# RECONSTRUCTING A TWO-COLOR SCENERY BY OBSERVING IT ALONG A SIMPLE RANDOM WALK PATH

By Heinrich Matzinger

*Universität Bielefeld and Georgia Tech*

Let $\{\xi(n)\}_{n\in\mathbb{Z}}$ be a two-color random scenery, that is, a random coloring of $\mathbb{Z}$ in two colors, such that the $\xi(i)$'s are i.i.d. Bernoulli variables with parameter $\frac{1}{2}$. Let $\{S(n)\}_{n\in\mathbb{N}}$ be a symmetric random walk starting at 0. Our main result shows that a.s., $\xi\circ S$ (the composition of $\xi$ and $S$) determines $\xi$ up to translation and reflection. In other words, by observing the scenery $\xi$ along the random walk path $S$, we can a.s. reconstruct $\xi$ up to translation and reflection. This result gives a positive answer to the question of H. Kesten of whether one can a.s. detect a single defect in almost every two-color random scenery by observing it only along a random walk path.

**1. Introduction.** A scenery is defined to be a function from $\mathbb{Z}$ to $\{0,1\}$. Let $\xi$ and $\tilde{\xi}$ be two sceneries. We say that $\xi$ and $\tilde{\xi}$ are equivalent iff there exist $a\in\mathbb{Z}$ and $b\in\{-1,1\}$ such that for all $x\in\mathbb{Z}$ we have $\xi(x)=\tilde{\xi}(a+bx)$. In this case we write $\xi\approx\tilde{\xi}$. In other words, two sceneries are equivalent iff they can be obtained from each other by a shift or a reflection. In everything that follows $\{S(k)\}_{k\geq 0}$ will be a simple random walk on $\mathbb{Z}$ starting at the origin. We will denote by $\chi\in\{0,1\}^{\mathbb{N}}$ the color record obtained by observing the scenery $\xi$ along the path of the random walk $\{S(k)\}_{k\geq 0}$:

$$\chi:=(\xi(S(0)),\xi(S(1)),\xi(S(2)),\dots),$$

that is, $\chi(k):=\xi(S(k))$ for all $k\in\mathbb{N}$. We examine the following question: given an unknown scenery $\xi$, can we "reconstruct" $\xi$ if we can only observe $\chi$? Thus, does one path realization of the process $\{\chi(k)\}_{k\geq 0}$ uniquely determine $\xi$? The answer in those general terms is "no." However, under appropriate restrictions, the answer will become "yes." This is the main result of this paper. Let us explain these restrictions: First, if $\xi$ and $\tilde{\xi}$ are equivalent, we can in general not distinguish whether the observations come from $\xi$ or from $\tilde{\xi}$. Thus, we can only reconstruct $\xi$ up to equivalence modulo $\approx$. Second, it is clear that the reconstruction will in the best case work only almost surely. If the random walk $\{S(k)\}_{k\geq 0}$ decides to walk only to the left (which it could do with probability zero), then we obtain no information about the right-hand side of the scenery $\xi$ and thus are not able to reconstruct the scenery $\xi$. Eventually, Lindenstrauss in [17] exhibits sceneries which one cannot reconstruct. Thus, not all sceneries can be reconstructed. However, we prove that many "typical" sceneries can be reconstructed up to equivalence and almost surely. For this we take the scenery $\xi$ to be the outcome of a random process which is independent of $\{S(k)\}_{k\geq 0}$ such that the $\xi(k)$'s are i.i.d. Bernoulli with parameter $\frac{1}{2}$. We use the following notation: we write $\xi$ for the (random) scenery: $\xi : k \mapsto \xi(k), \mathbb{Z} \to \{0,1\}$. Our main result states that, given only the observation $\chi$, almost every scenery $\xi$ can be reconstructed a.s. up to equivalence. Let us state our main theorem:

THEOREM 1. *Let $\{S(k)\}_{k\geq 0}$ and $\{\xi(k)\}_{k\in\mathbb{Z}}$ be two processes independent of each other such that $\{S(k)\}_{k\geq 0}$ is a simple random walk starting at the origin and such that the $\xi(k)$'s are i.i.d. Bernoulli variables with parameter $\frac{1}{2}$. Then a.s. $\chi$ determines $\xi$ up to equivalence. In other words, there exists a measurable function $\mathcal{A} : \{0,1\}^{\mathbb{N}} \to \{0,1\}^{\mathbb{Z}}$ such that $P(\mathcal{A}(\chi) \approx \xi) = 1$. ("Measurable" means measurable with respect to the $\sigma$-algebras induced by the canonical coordinates on $\{0,1\}^{\mathbb{N}}$ and on $\{0,1\}^{\mathbb{Z}}$.)*

We will prove the above theorem by explicitly describing how to reconstruct $\xi$ from $\chi$. Hence, our approach is constructive. We explicitly give a construction which produces a (random) scenery $\bar{\xi} : \mathbb{Z} \to \{0,1\}$ when applied to the observations $\chi$. The constructed scenery $\bar{\xi}$ is shown to be a.s. equivalent to $\xi$. In this way $\mathcal{A}$ gets defined: $\mathcal{A}(\chi) := \bar{\xi}$.

Let us now make a few historical comments. This paper was motivated by Kesten's question to me of whether one can a.s. distinguish a single defect in almost any two-color scenery. Let us explain what the scenery distinguishing problem is. Let $\xi, \eta : \mathbb{Z} \to \{0,1\}$ and let $\{S(k)\}_{k\in\mathbb{N}}$ be a symmetric random walk on $\mathbb{Z}$. Let the process $\{\chi(k)\}_{k\in\mathbb{N}}$ be equal to either $\{\xi(S(k))\}_{k\in\mathbb{N}}$ or $\{\eta(S(k))\}_{k\in\mathbb{N}}$. Is it possible by observing only one path realization of $\{\chi(k)\}_{k\in\mathbb{N}}$ to say to which one of the two $\{\xi(S(k))\}_{k\in\mathbb{N}}$ or $\{\eta(S(k))\}_{k\in\mathbb{N}}$, $\{\chi(k)\}_{k\in\mathbb{N}}$ is equal to? (We assume that we know $\xi$ and $\eta$.) If yes, we say that it is possible to distinguish between the sceneries $\xi$ and $\eta$ by observing them along a path of $\{S(k)\}_{k\in\mathbb{N}}$. Otherwise, when it is not possible to figure out almost surely by observing $\{\chi(k)\}_{k\in\mathbb{N}}$ alone whether $\{\chi(k)\}_{k\in\mathbb{N}}$ is generated on $\xi$ or on $\eta$, we say that $\xi$ and $\eta$ are indistinguishable. The problem of distinguishing two sceneries was raised independently by Benjamini and by den Hollander and Keane. The motivation came from problems in ergodic theory, such as the $T, T^{-1}$problem (see [10]) and from the study of various aspects of $\{\xi(S(k))\}_{n\in\mathbb{N}}$, where $\{\xi(k)\}_{k\in\mathbb{Z}}$ is random. (See [3, 11, 14].) Benjamini and Kesten showed in [1] that one can distinguish almost any two random sceneries even when the random walk is in $\mathbb{Z}^2$. (They assumed the sceneries to be random themselves, so that the $\xi(k)$'s and the $\eta(n)$'s are i.i.d. Bernoulli.) Kesten in [12] proved that when the random sceneries are i.i.d. and have four colors, that is, $\xi$ and $\eta:\mathbb{Z}\to\{0,1,2,3\}$, and differ only in one point, they can be a.s. distinguished. He asked whether this result might still hold with fewer colors. The main result of this paper directly implies that one can distinguish single defects in almost any scenery. In [21], we proved for the three-color case that one can a.s. reconstruct almost every three-color scenery. We also established that this implies that one can distinguish single defects for almost all three-color sceneries. In the two-color case, that is, in the case we consider in this paper, the same thing is true. This means that our result for scenery reconstruction implies that one can distinguish single defects in almost all sceneries. We state the following corollary to our main result without giving a proof. (The proof that our main result implies the following corollary is very similar to the one given in [21] for the three-color case.)

COROLLARY 2. *Let $\mathbb{B}$ designate the set of all two-color sceneries.* $\mathbb{B}=\{\xi:\mathbb{Z}\to\{0,1\}\}=\{0,1\}^{\mathbb{Z}}$. *Let $(\mathbb{B},\sigma(\mathbb{B}))$ denote the measurable space, where $\sigma(\mathbb{B})$ is the $\sigma$-algebra induced by the canonical coordinates on $\mathbb{B}$. Let $P$ denote the probability measure on $(\mathbb{B},\sigma(\mathbb{B}))$ obtained by assuming that the $\xi(i)$'s are i.i.d. Bernoulli variables with parameter $\frac{1}{2}$. Then there exists a $\sigma(\mathbb{B})$-measurable set $\mathbb{S}$, such that $P(\mathbb{S})=1$ and such that for every scenery $\xi\in\mathbb{S}$ and every scenery $\eta$ which is equal to $\xi$ everywhere except in one point, we have that $\xi$ and $\eta$ are distinguishable.*

The above corollary says that there are many sceneries which one can distinguish or, in other words, that sceneries which are typical in a certain sense can be distinguished. However, the above result becomes false if one tries to extend it to all pairs of sceneries which are not equivalent. Recently, Lindenstrauss [17] exhibited a nondenumerable set of pairs of nonequivalent sceneries on $\mathbb{Z}$ which he proved to be indistinguishable. Before that, Howard proved in [[7]–[9]] that any two periodical sceneries of $\mathbb{Z}$ which are not equivalent modulo translation and reflection are distinguishable and that one can a.s. distinguish single defects in periodical sceneries. Kesten asked in [13] whether this result would still hold when the random walk would be allowed to jump. He also asked what would happen in the two-dimensional case. Löwe and Matzinger in [18] have been able to prove that one can a.s. reconstruct almost every scenery up to equivalence in two dimensions, provided the scenery has a lot of colors. However, the problem of the reconstruction of two-color sceneries in $\mathbb{Z}$ seen along the random walk path of a recurrent random walk which is allowed to jump remains open. In our opinion, this is a central open problem at present. Eventually we should also mention that the two-color scenery reconstruction problem for a scenery which is i.i.d. is equivalent to the following problem: let $\{R(k)\}_{k\in\mathbb{Z}}$ and $\{S(k)\}_{k\geq 0}$ be two independent simple random walks on $\mathbb{Z}$ both starting at the origin and living on the same probability space. (Here we mean that $\{R(k)\}_{k\geq 0}$ and $\{R(-k)\}_{k\geq 0}$ are two independent simple random walks both starting at the origin.) Does one path realization of the iterated random walk $\{R(S(k))\}_{k\geq 0}$ uniquely determine the path of $\{R(k)\}_{k\in\mathbb{Z}}$ up to shift and reflection around the origin? If one takes the representation of the scenery $\xi$ as a nearest-neighbor walk (which we will define later) for $\{R(k)\}_{k\in\mathbb{Z}}$, then it becomes immediately clear that the two problems are equivalent. We leave it to the reader to check the details. So the main result of this paper is equivalent to the following result for iterated nearest-neighbor walks: one path realization of the iterated random walk $\{R(S(k))\}_{k\geq 0}$ a.s. uniquely determines the path of $\{R(k)\}_{k\in\mathbb{Z}}$ up to shift and reflection around the origin. This is a discrete analog of the result of Burdzy [2] concerning the path of iterated Brownian motion.

**2. Reconstructing a finite piece of the scenery $\boldsymbol{\xi}$.** To explain a key idea, we first present a solution to a simplified but somewhat unrealistic case.

2.1. *Simplified example.* Assume for a moment that the scenery $\xi$ is nonrandom, and instead of being a two-color scenery, is a four-color scenery, that is, $\xi:\mathbb{Z}\to\{0,1,2,3\}$. Let us imagine furthermore, that there are two integers $x,y$ such that $\xi(x)=2$ and $\xi(y)=3$, but outside $x$ and $y$ the scenery has everywhere color 0 or 1, [i.e., for all $z\in\mathbb{Z}$ with $z\neq x,y$ we have that $\xi(z)\in\{0,1\}$]. The simple random walk $\{S(k)\}_{k\geq 0}$ can go with each step one unit to the right or one unit to the left. This implies that the shortest possible time for the random walk $\{S(k)\}_{k\geq 0}$ to go from the point $x$ to the point $y$ is $|x-y|$. When the random walk $\{S(k)\}_{k\geq 0}$ goes in shortest possible time from $x$ to $y$ it goes in a straight way, which means that between the time it is at $x$ and until it reaches $y$ it only moves in one direction. During that time, the random walk $\{S(k)\}_{k\geq 0}$ reveals the portion of $\xi$ lying between $x$ and $y$. If between time $t_1$ and $t_2$ the random walk goes in a straight way from $x$ to $y$ [i.e., if $|t_1-t_2|=|x-y|$ and $S(t_1)=x, S(t_2)=y$], then the word $\chi(t_1),\chi(t_1+1),\dots,\chi(t_2)$ is a copy of the scenery $\xi$ restricted to the interval $[\min\{x,y\},\max\{x,y\}]$. In this case, the word $\chi(t_1),\chi(t_1+1),\dots,\chi(t_2)$ is equal to the word $\xi(x),\xi(x+u),\xi(x+2u),\dots,\xi(y)$, where $u:=(y-x)/|y-x|$. Since the random walk $\{S(k)\}_{k\geq 0}$ is recurrent, it a.s. goes at least once, in the shortest possible way, from the point $x$ to the point $y$. Because we are given infinitely many observations we can a.s. figure out what the distance between $x$ and $y$ is: the distance between $x$ and $y$ is the shortest time lapse that a "3" will ever appear in the observations $\chi$ after a "2." When, on the other hand, a "3" appears in the observations $\chi$ in shortest possible time after a "2," then between the time we see that "2" and until we see the next "3," we observe a copy of $\xi(x), \xi(x+u), \xi(x+2u), \dots, \xi(y)$ in the observations $\chi$. This fact allows us to reconstruct the finite piece $\xi(x), \xi(x+u), \xi(x+2u), \dots, \xi(y)$ of the scenery. Choose any couple of integers $t_1, t_2$ with $t_2 > t_1$, minimizing $|t_2 - t_1|$ under the condition that $\chi(t_1) = 2$ and $\chi(t_2) = 3$. Almost surely then $\chi(t_1), \chi(t_1+1), \dots, \chi(t_2)$ is equal to $\xi(x), \xi(x+u), \xi(x+2u), \dots, \xi(y)$.

A NUMERICAL EXAMPLE. Let the scenery $\xi$ be such that $\xi(-2) = 0$, $\xi(-1) = 2$, $\xi(0) = 0$, $\xi(1) = 1$, $\xi(2) = 1$, $\xi(3) = 3$, $\xi(4) = 0$. Assume furthermore that the scenery $\xi$ has a 2 and a 3 nowhere else than in the points $-1$ and 3. Imagine that $\chi$, the observations given to us, would start as follows:

$$\chi = (0, 2, 0, 1, 0, 1, 1, 3, 0, 3, 1, 1, 1, 1, 0, 2, 0, 1, 1, 3, \dots).$$

By looking at all of $\chi$ we would see that the shortest time a 3 occurs after a 2 in the observations is 4. In the first observations given above there is, however, already a 3 only four time units after a 2. The binary word appearing in that place, between the 2 and the 3, is 011. We deduce from this that between the place of the 2 and the 3 the scenery must look like: 011.

In reality the scenery we want to reconstruct has two colors only. So, instead of the 2 and the 3 in the example above we will use a special pattern in the observations which will tell us when the random walk is back at the same spot. One possibility (although not yet the one we will eventually use) would be to use binary words of the form: 001100 and 110011. It is easy to verify that the only possibility for the word 001100, respectively, 110011, to appear in the observations is when the same word 001100, respectively, 110011, occurs in the scenery and the random walk reads it. So, imagine (to give another pedagogical example of a simplified case) the scenery would be such that in a place $x$ there occurs the word 001100, and in the place $y$ there occurs the word 110011, but these two words occur in no other place in the scenery. These words can then be used as markers: In order to reconstruct the piece of the scenery $\xi$ included between $x$ and $y$ we could proceed as follows: take in the observations the place where the word 110011 occurs in shortest time after the word 001100. In that place in the observations we see a copy of the piece of the scenery $\xi$ included between $x$ and $y$. The reason why the very last simplified example is not realistic is the following: we take the scenery to be the outcome of a random process itself where the $\xi(k)$'s are i.i.d. variables themselves. Thus any word will occur infinitely often in the scenery $\xi$. However, if, for example, the markers in the scenery occur far away from each other, then we can still use the above described reconstruction strategy: The random walk will then be very likely to first cross from $x$ to $y$ in a straight way before meeting another marker and creating some confusion. In the next section we explain how to construct the markers which we are eventually going to use.

2.2. *Representation of the scenery $\xi$ as a nearest-neighbor walk.* The scenery reconstruction problem contains two main ingredients: A random walk $\{S(k)\}_{k\in\mathbb{N}}$ and a "random environment," that is, the scenery $\xi$. The key idea in this paper is to view the random environment itself as a nearest-neighbor walk. In this section we explain how to do this, by defining "the representation of the scenery $\xi$ as a nearest-neighbor walk." We need the following definitions: Let $D$ be an integer interval, that is, the intersection between a real interval and the integer numbers $\mathbb{Z}$. We call a function $T: D \to \mathbb{Z}$ a nearest-neighbor walk, iff for each $t_1, t_2 \in D$ with $|t_1 - t_2| = 1$, we have that $|T(t_1) - T(t_2)| = 1$. In what follows, we will write $S$ for the path of the process $\{S(k)\}_{k\geq 0}$, that is, for $S: k \mapsto S(k), \mathbb{N} \to \mathbb{Z}$. Let $\varphi: \mathbb{Z} \to \{0,1\}$ be one of the two 4-periodic sceneries with period 0011 and $\varphi(0) = \varphi(1)$. Such a scenery $\varphi$ has a very particular property: for every point in the scenery $\phi$, one neighboring point has color 0, while the other one has color 1. This implies that for any color record $\phi$ there exists one and only one nearest-neighbor walk $T$ generating $\phi$ on the scenery $\varphi$ once we know where $T$ starts. We can use this fact to represent a color record as a nearest-neighbor walk: the nearest-neighbor walk representing a sequence of colors is simply defined to be the only nearest-neighbor walk generating the color sequence on $\varphi$ and starting at a given point, in general the origin. (For this to work the starting point must have the right color.)

A NUMERICAL EXAMPLE. Let $\phi = (01011000101010100\dots)$ be a color record we want to represent as a nearest-neighbor walk. Let $\varphi: \mathbb{Z} \to \{0,1\}$ be the 4-periodic scenery:

| $\varphi(k)$ | … | 0 | 0 | 1 | 1 | 0 | 0 | 1 | 1 | 0 | 0 | 1 | 1 | … |
|---|---|---|---|---|---|---|---|---|---|---|---|---|---|---|
| $k$ | … | −4 | −3 | −2 | −1 | 0 | 1 | 2 | 3 | 4 | 5 | 6 | 7 | … |

Define the nearest-neighbor walk representing $\phi$ to be the only nearest-neighbor walk $T: \mathbb{N} \to \mathbb{Z}$ starting at the origin and generating the sequence $\phi$ on $\varphi$, that is, such that $\varphi \circ T = \phi$. In this example we get

| $T(t)$ | 0 | −1 | 0 | −1 | −2 | −3 | −4 | −3 | −2 | −3 | −2 | … |
|---|---|---|---|---|---|---|---|---|---|---|---|---|
| $t$ | 0 | 1 | 2 | 3 | 4 | 5 | 6 | 7 | 8 | 9 | 10 | … |

The scenery $\xi$ we want to represent as a nearest-neighbor walk is, however, a doubly infinite sequence. We will thus take the sequence $\xi(0), \xi(1), \xi(2), \xi(3), \dots$ first and define with it the portion of the path of the nearest-neighbor walk in positive time. Then we take $\xi(0), \xi(-1), \xi(-2), \xi(-3), \dots$, and this defines the part of the nearest-neighbor walk in negative time.

An example. Let $\xi : \mathbb{Z} \to \{0, 1\}$ be a scenery with the following values close to the origin:

| $\xi(k)$ | $\dots$ | 1 | 0 | 1 | 0 | 0 | 0 | 1 | 1 | 1 | 0 | 0 | 1 | $\dots$ |
|---|---|---|---|---|---|---|---|---|---|---|---|---|---|---|
| $k$ | $\dots$ | $-4$ | $-3$ | $-2$ | $-1$ | 0 | 1 | 2 | 3 | 4 | 5 | 6 | 7 | $\dots$ |

Designate by $R$ the nearest-neighbor walk representing $\xi$. Then the part of $\xi$ to the right of the origin defines the path of $R$ which lies in positive time. In this example, $(00111001\dots)$ is responsible for this part of $R$. We get

| $R(t)$ | 0 | 1 | 2 | 3 | 2 | 1 | 0 | $-1$ | $\dots$ |
|---|---|---|---|---|---|---|---|---|---|
| $t$ | 0 | 1 | 2 | 3 | 4 | 5 | 6 | 7 | $\dots$ |

In the same way, the part of $\xi$ which lies left of the origin is responsible for the restriction of $R$ to the negative integers. In our example, $(\dots 1010)$ defines that part of $R$. We get

| $R(t)$ | $\dots$ | 2 | 1 | 2 | 1 | 0 |
|---|---|---|---|---|---|---|
| $t$ | $\dots$ | $-4$ | $-3$ | $-2$ | $-1$ | 0 |

We are ready to define $R$ formally:

Definition 3. Let $\varphi : \mathbb{Z} \to \{0, 1\}$ designate the following 4-periodic (random) scenery:

1. When $\xi(0) = 0$, we set $(\varphi(0), \varphi(1), \varphi(2), \varphi(3)) = (0, 0, 1, 1)$.
2. When $\xi(0) = 1$, we set $(\varphi(0), \varphi(1), \varphi(2), \varphi(3)) = (1, 1, 0, 0)$.

The nearest-neighbor walk $R : \mathbb{Z} \to \mathbb{Z}$ representing the scenery $\xi$ is defined to be the only (random) nearest-neighbor walk $R$ such that $R(0) = 0$ and $\varphi \circ R = \xi$, that is, $\varphi(R(k)) = \xi(k)$ for all $k \in \mathbb{Z}$.

It is easy to verify that when the $\xi(k)$'s are i.i.d. Bernoulli variables with $P(\xi(0) = 0) = P(\xi(0) = 1) = \frac{1}{2}$, then $\{R(k)\}_{k \in \mathbb{Z}}$ as well as $\{R(-k)\}_{k \in \mathbb{Z}}$ are two independent symmetric random walks starting at the origin.

In Figure 1, we illustrate the above numerical example by showing a portion of the graph of $R$. For this we take

$$
\begin{aligned}
&(\xi(0), \xi(1), \xi(2), \xi(3), \xi(4), \dots) \\
&\quad = (00111001011000110010000100110010010110010011100 1\dots).
\end{aligned}
$$

In Figure 1, the label $k$ designates the point $(R(k), k)$.

Next we need a few definitions.

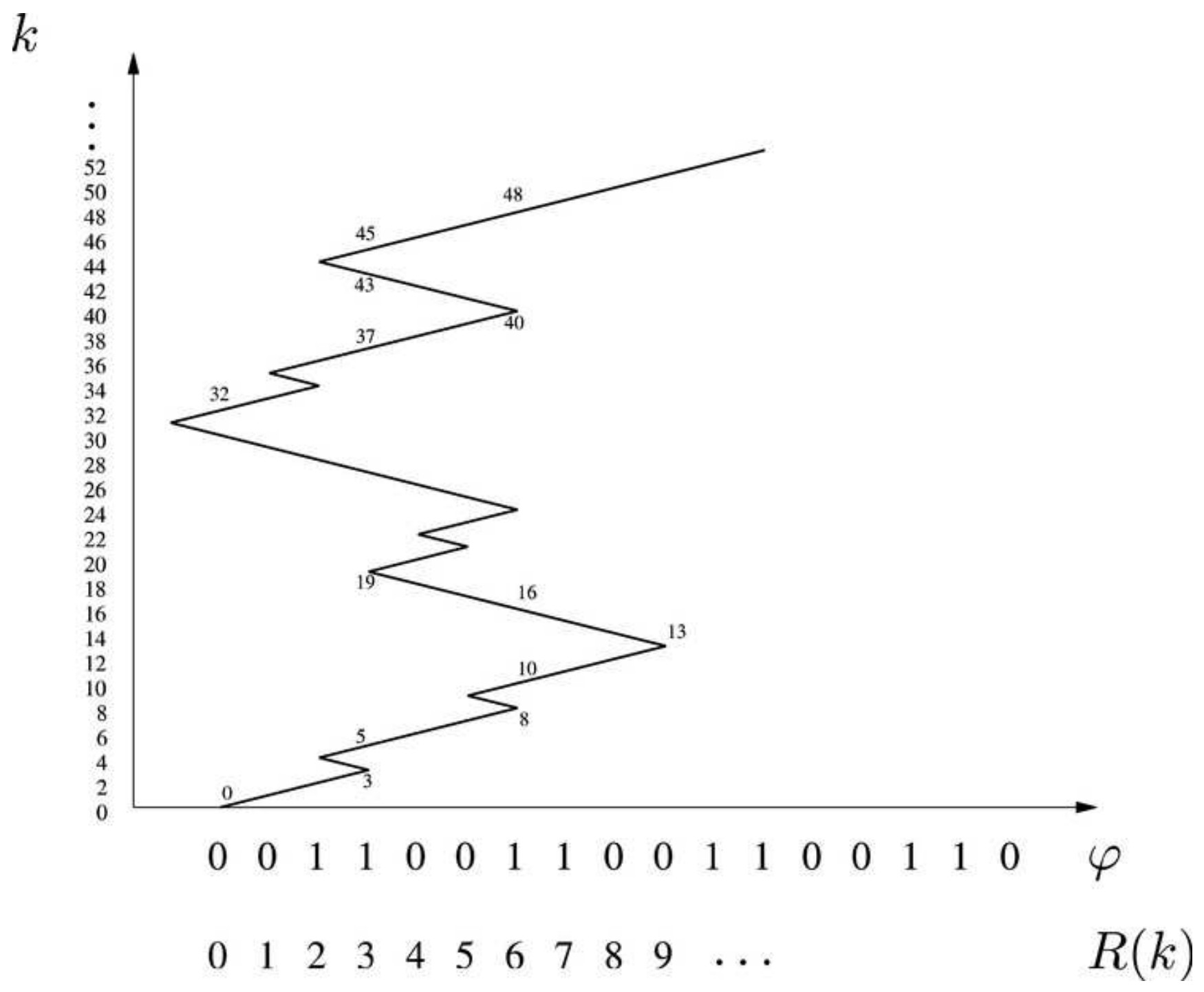

FIG. 1.

DEFINITION 4. Let $T : D \to \mathbb{Z}$ be a nearest-neighbor walk. Let $t_1$, $t_2 \in D$ and $x_1, x_2 \in \mathbb{Z}$, $x_1 \neq x_2$. We call $(t_1, t_2)$ a crossing by $T$ of $(x_1, x_2)$ iff $(T(t_1), T(t_2)) = (x_1, x_2)$ and for all integer $t$ strictly between $t_1$ and $t_2$, $T(t)$ is strictly between $x_1$ and $x_2$. If $t_2 > t_1$ we say that the crossing $(t_1, t_2)$ is "positive," otherwise we say that it is "negative." If $|t_1 - t_2| = |x_1 - x_2|$, we say that the crossing $(t_1,\, t_2)$ is straight.

Let $(t_3,\, t_4)$ be a crossing by $T$ of $(x_3, x_4)$. Then, we say that $(t_3,\, t_4)$ is the *first crossing* by $T$ of $(x_3, x_4)$ during $(t_1, t_2)$ iff $t_3,\, t_4 \in [\min\{t_1, t_2\}, \max\{t_1, t_2\}]$ and $(t_3,\, t_4)$ is the crossing by $T$ of $(x_3, x_4)$ which lies in $[\min\{t_1, t_2\}, \max\{t_1, t_2\}]$ (i.e., $t_3,\, t_4 \in [\min\{t_1, t_2\}, \max\{t_1, t_2\}]$) and is *closest to* $t_1$.

Let $(t_1, t_2)$ and $(s_1, s_2)$ be two crossings by nearest-neighbor walk $T$ of $(x_1, x_2)$. Then either the intervals:

$$]\min\{t_1, t_2\}, \max\{t_1, t_2\}[ \quad \text{and} \quad ]\min\{s_1, s_2\}, \max\{s_1, s_2\}[$$

are disjoint, or $(t_1, t_2) = (s_1, s_2)$ holds. Thus, we can numerate the crossings by $T$ of $(x_1, x_2)$ in increasing order of appearance. Thus the above definition of "first crossing by $T$ of $(x_3, x_4)$ during another crossing" makes sense.

In the numerical example of Figure 1, we see that between time 0 and time 3 the nearest-neighbor walk $R$ crosses from the point 0 to the point 3 in a straight way. In other words, $(0,3)$ is a straight crossing by $R$ of $(0,3)$. Furthermore, $R$ during the time interval $(0,13)$ crosses the interval $(0,9)$. Thus, $(0,13)$ is a crossing by $R$ of $(0,9)$. Because $(0,3) \in (0,13)$ we have that the crossing $(0,3)$ happens during the crossing $(0,13)$. Clearly, $(0,3)$ is the first crossing by $R$ of $(0,3)$ during the crossing $(0,13)$. (In the above example it is also the only one.) The crossing $(0,13)$, unlike $(0,3)$, is not a straight one. $(32,51)$ is a crossing by $R$ of $(0,9)$. This is the second crossing by $R$ of $(0,9)$ after time 0. During the crossing $(32,51)$ there are two crossings by $R$ of the $(3,6)$. These are $(37,40)$ and $(45,48)$.

2.3. *Localization test.* In this section, we construct a test to determine at what times the random walk is back at the same location. Combined with the idea of "going in shortest time from $x$ to $y$," we have the main ingredients for the reconstruction of a finite piece of the scenery $\xi$. If we have such a test, we can recognize when the random walk is back at a location $x$ and at which times it is back at locations $x$ and $y$. We then take a time interval where the random walk visits $y$ in shortest possible time after visiting $x$.

This "localization test" is based on the representation $R$ of the scenery $\xi$ as a nearest-neighbor walk. Recall that $R$ is not observable. The composition of two nearest-neighbor walks is again a nearest-neighbor walk. Thus, the composition $R \circ S : k \mapsto R(S(k)), \mathbb{N} \to \mathbb{Z}$ is a nearest-neighbor walk. However, every nearest-neighbor walk $T : \mathbb{N} \to \mathbb{Z}$ is uniquely determined by $\varphi \circ T$. In the following we set

$$T := R \circ S.$$

We get

$$\varphi \circ T = (\varphi \circ R) \circ S = \xi \circ S = \chi,$$

that is, $T$ generates the color record $\chi$ on the scenery $\varphi$. Furthermore, $T(0) = 0$. Thus $T$ is uniquely determined by the observations $\chi$. Hence $T$ is observable. Thus, although $R$ and $S$ are both not known, their composition $R \circ S$ is observable. We are using the nearest-neighbor walk $R \circ S$ to determine when $S$ is back at the same place.

To illustrate the *main idea* of the localization test (and maybe of this paper) we view the random walk $S$ on the graph $k \mapsto (R(k), k)$ geometrically in two dimensions. This defines a movement in two dimensions:

$$t \mapsto (R(S(t)), S(t)).$$

By projecting this movement along the $y$-axis on the $x$-axis we get the known one-dimensional nearest-neighbor walk $T$. Imagine that the path of $R$ is given; then $t \mapsto (R(S(t)), S(t))$ can be viewed as a one-dimensional random walk moving in $\mathbb{R}^2$ on the graph of $R$.

Figure 2 illustrates this situation. The graph of $R$ is drawn as a dotted line, as it is not observable. The hand-drawn lines with arrows represent the movement of the random walk $S$ on the graph of $G$. This is the movement $t \mapsto (R(S(t), S(t)))$, which is also not observable. However, the projection of this movement onto the horizontal line gives the observable nearest-neighbor walk $R \circ S$, which is observable.

Let $\Delta S(k) := S(k+1) - S(k)$. In the example of Figure 2 we have that

$$(\Delta S(0), \Delta S(1), \Delta S(2), \dots)$$
$$= (+1, +1, +1, +1, +1, +1, +1,$$
$$-1, +1, +1, +1, +1, +1, -1, +1, +1, +1, \dots)$$

and $R$ takes on the same values as in Figure 1.

Imagine that the dotted line representing the graph of $R$ is made out of invisible glass. The random walk $S$ moves invisibly on that glass line, but its projection onto the $x$-axis is visible. Seeing only this projection, we want to determine when $S$ has returned to the same place. $S$ has returned exactly when the two-dimensional movement $t \mapsto (R(S(t)), S(t))$ has returned to the same place: $S(s) = S(t)$ iff $(R(S(s)), S(s)) = (R(S(t)), S(t))$. Viewing $R$ as fixed, this means that $S$ is back at the same place exactly when the random walk $S$ on the graph of $R$ has come back to the same place. As shown below, we can statistically determine this with high precision by counting the number of straight crossings of $R \circ S$ and their location. Let us illustrate the idea with Figure 3.

In Figure 3, we show two finite portions of the movement of the random walk $S$ on the graph of $R$. The first one is designated by the letter $a$ while the second one is designated by the letter $b$. In this example $a$ corresponds to the random walk $S$ making the following first steps:

$$(\Delta S(0), \Delta S(1), \Delta S(2), \dots)$$
$$= (+1, +1, +1, +1, +1, +1, +1,$$
$$-1, +1, +1, +1, +1, +1, -1, +1, +1, +1, \dots).$$

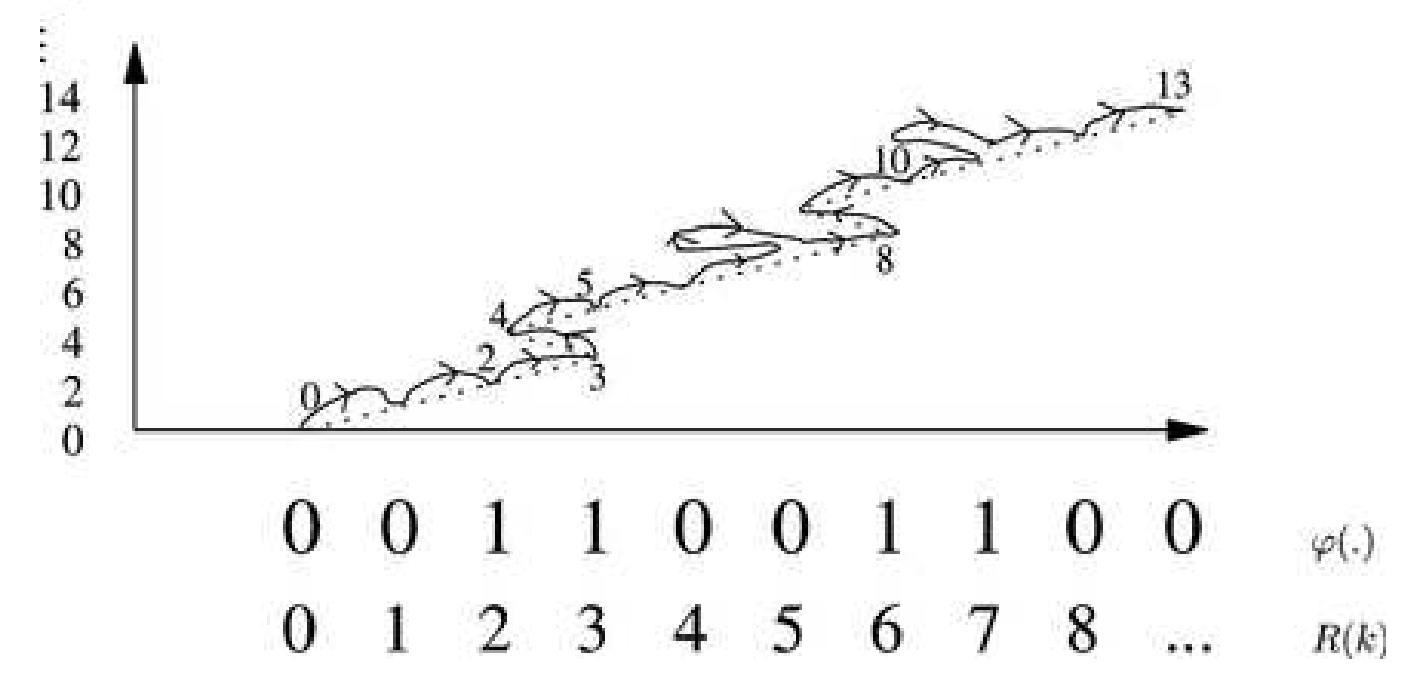

Fig. 2.

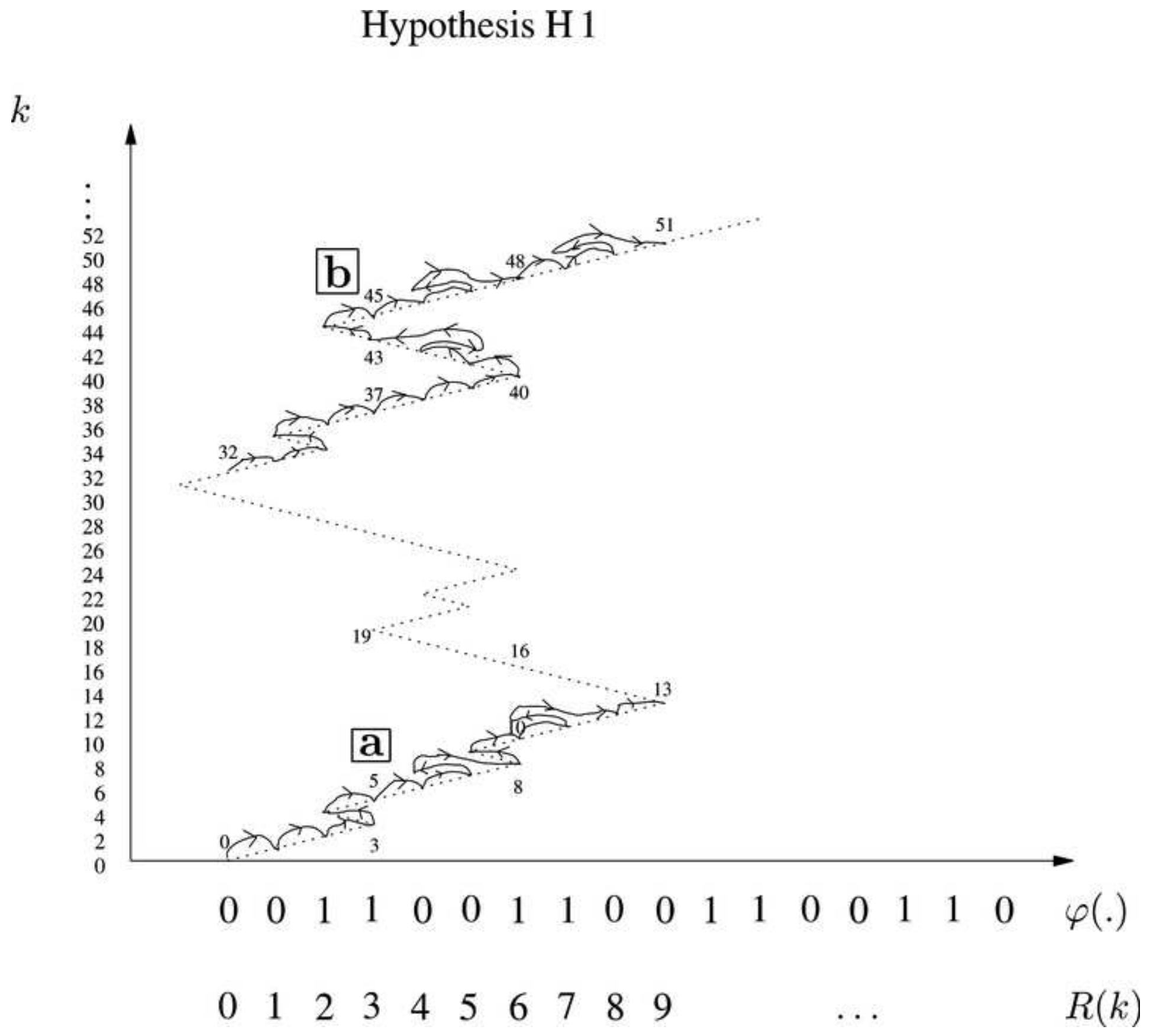

FIG. 3.

Part $b$ starts at time $t_b$ such that $S(t_b) = 32$. Then the random walk $S$ makes the following steps:

$$(\Delta S(t_b), \Delta S(t_b+1), \Delta S(t_b+2), \dots)$$
$$= (+1, +1, +1, +1, +1, +1, +1, +1, +1, +1,$$
$$-1, +1, +1, +1, +1, +1, +1, -1, +1, +1, +1, +1, -1, +1, +1 \dots).$$

The random walk $S$ from time $t_b$ until time $t_b + 25$ performs a crossing of the interval $(32, 51)$. This means that at time $t_b$ the random walk $S$ is at the point 32 and at time $t_b + 25$ it is at the point 51, but strictly in between the time $t_b$ until time $t_b + 25$ the random walk $S$ does not visit the points 32 or 51. In Figure 3 if we project the movement $b$ (of the random walk $S$ on the graph of $R$) onto the horizontal line, we get the movement of the nearest-neighbor walk $R \circ S$ during the time interval from time $t_b$ until time $t_b + 25$. This is a crossing as well: during that time $R \circ S$ crosses from the point 0 to the point 9; that is, it crosses the interval $(0,9)$. During that time $S$ on the graph of $R$ crosses a portion of the graph of $R$ which corresponds itself to a crossing by $R$. As a matter of fact, between time 32 and time 51 the nearest-neighbor walk $R$ crosses the interval $(0,9)$. Following our convention we say that $(32,51)$ is a crossing by the nearest-neighbor walk $R$ of the interval $(0,9)$. In part $a$ we see the following: $(0,17)$ is a crossing by $S$ of $(0,13)$. On the other hand, $(0,13)$ is a crossing by $R$ of $(0,9)$. Eventually, $(0,17)$ is a crossing by $R\circ S$ of $(0,9)$.

The example of Figure 3 illustrates one of the three main combinatorial facts used in this paper: the composition $T=R\circ S$ performs a crossing iff during that time $S$ performs a crossing of a crossing of $R$. Let us formulate this as a lemma:

LEMMA 5. *Let $0<t_1<t_2$. $(t_1,t_2)$ is a crossing by $T$ of the interval $(x_1,x_2)$ iff there exist $k_1,k_2\in\mathbb{Z}$ such that $(t_1,t_2)$ is a crossing by $S$ of $(k_1,k_2)$, and $(k_1,k_2)$ is a crossing by $R$ of $(x_1,x_2)$.*

Let us study next the example of Figure 3 more: during time $(14,17)$, $S$ performs a straight crossing of the interval $(10,13)$. Furthermore, $(10,13)$ represents itself a straight crossing by $R$ of the interval $(6,9)$. This leads to, that $R\circ S$ performs during the time interval $(14,17)$ a straight crossing of the interval $(6,9)$. On the other hand, during time $(t_b,t_b+4)$ $S$ performs a straight crossing of the interval $(32,37)$. However, $(32,37)$ is a crossing by $R$, but not a straight one. It follows that $(t_b,t_b+4)$ is a crossing by $R\circ S$, but not a straight one.

The rule is: on a crossing by $R$ which is not straight it is impossible to get a crossing by $R\circ S$ which is straight. This is the second main combinatorial fact:

LEMMA 6. *Let $0<t_1<t_2$. Then $(t_1,t_2)$ is a straight crossing by $T$ of the interval $(x_1,x_2)$ iff there exists $k_1,k_2\in\mathbb{Z}$ such that $(t_1,t_2)$ is a straight crossing by $S$ of $(k_1,k_2)$ and $(k_1,k_2)$ is a straight crossing by $R$ of $(x_1,x_2)$.*

Looking further at Figure 3, we see that in portion $b$ of the path of $S$ on the graph of $R$ we have: during the crossing $(32,51)$ the first crossing by $R$ of $(3,6)$ is $(37,40)$ and the last one is $(45,48)$. The first crossing by $S$ of $(37,40)$ during $(t_b,t_b+51)$ is $(t_b+5,t_b+8)$. The first crossing during $(t_a,t_a+25)$ by $R\circ S$ of $(3,6)$ is also $(t_b+5,t_b+8)$. Thus, the first crossing during $(t_a,t_a+25)$ by $R\circ S$ of $(3,6)$ happens when during $(t_a,t_a+25)$ $S$ crosses for the first time the first crossing by $R$ of $(3,6)$.

We see that a first crossing by $R\circ S$ corresponds to a first crossing by $S$ of a first crossing by $R$. This yields our third combinatorial fact:

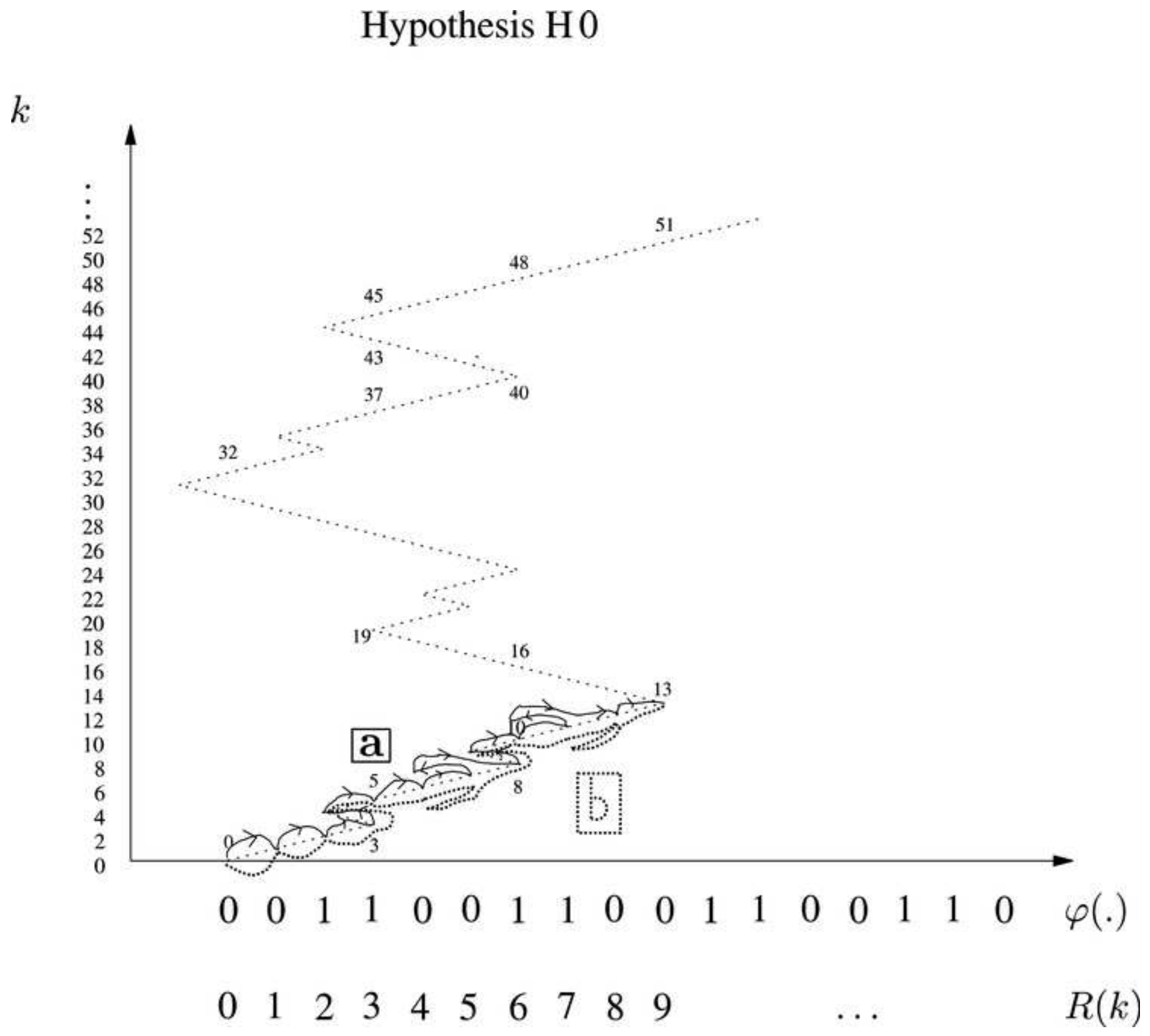

FIG. 4.

LEMMA 7. *Let $0 < t_1 < t_2 < t_3 < t_4$ and $0 < x_1 < x_2 < x_3 < x_4$. Furthermore, let $(t_1, t_4)$ be a crossing by $R \circ S$ of $(x_1, x_4)$. Then $(t_2, t_3)$ is the first crossing during $(t_1, t_4)$ of $(x_2, x_3)$ by $R \circ S$ iff it is the first crossing by $S$ during $(t_1, t_4)$ of $(k_2, k_3)$, where $(k_2, k_3)$ is the first crossing by $R$ of $(x_2, x_3)$ during $(k_1, k_4)$.*

To illustrate this, consider Figures 3 and 4.

In Figure 4, the portion $b$ of the path of the random walk $S$ is traced on the graph of $R$ as a thick dotted line. It is a crossing by $S$ of the crossing $(0, 13)$ by $R$. The projection onto the horizontal line of this movement is a crossing, too. In Figure 4, the crossing $b$, that is, $(t_b, t_b + 18)$, is a crossing by $S$ of $(0, 13)$. Furthermore, $(t_b, t_b + 18)$ is also a crossing by $R \circ S$ of $(0, 9)$.

Figure 4 is identical to Figure 3 except for the path $b$. In Figure 3, the crossings $a$ and $b$ by $S$ take place "in different locations," while in Figure 4 they take place "on the same locations." Given Figures 3 and 4, one can see if the crossings $a$ and $b$ by $S$ take place "in the same location" or not. However, from the input data of the reconstruction problem, only the projection down onto the horizontal of the path of $S$ on the graph of $R$ is observable: In both cases, we observe two crossings $a$ and $b$ by $R \circ S$ of the same interval $(0, 9)$. Based on the observation of those crossings only, we need to infer if "the crossings occur on the same location" as in Figure 4, or "on different locations" as in Figure 3.

More generally: Assume we observe two crossings $(t_{1a}, t_{2a})$ and $(t_{1b}, t_{2b})$ by $R \circ S$ of an interval $(0, 3n)$; this interval instead of any other interval is chosen for notational convenience. Because of Lemma 5, there exist $k_{1a}, k_{2a}$ such that $(k_{1a}, k_{2a})$ is a crossing by $R$ of $(0, 3n)$ while $(t_{1a}, t_{2a})$ is a crossing by $S$ of $(k_{1a}, k_{2a})$. Similarly, there exist $k_{1b}, k_{2b}$ such that $(k_{1b}, k_{2b})$ is a crossing by $R$ of $(0, 3n)$ while $(t_{1b}, t_{2b})$ is a crossing by $S$ of $(k_{1b}, k_{2b})$.

In Figure 3, $(t_{1a}, t_{2a}) = (0, 17)$, $(k_{1a}, k_{2a}) = (0, 13)$, $t_{2b} = t_{1b} + 25$, $(k_{1b}, k_{2b}) = (32, 51)$.

We develop a statistical test to determine if the two crossings $(t_{1a}, t_{2a})$ and $(t_{1b}, t_{2b})$ by $S$ occur "on the same place" or not. Its input data are two observed crossings $(t_{1a}, t_{2a})$ and $(t_{1b}, t_{2b})$ by $R \circ S$ of the same interval. We define the hypotheses of our test:

Hypothesis $H_0$. During the crossings $(t_{1a}, t_{2a})$ and $(t_{1b}, t_{2b})$ the random walk $S$ is on the same crossing of $R$. More precisely, $(S(t_{1a}), S(t_{2a})) = (S(t_{1b}), S(t_{2b}))$.

Hypothesis $H_1$. $(S(t_{1a}), S(t_{2a})) \neq (S(t_{1b}), S(t_{2b}))$.

If $H_0$ holds, then $S(t_{2a}) = S(t_{2b})$, that is, the random walk is back at the same place.

To determine if during two crossings by $R \circ S$ the random walk $S$ was at the same place we are going to count the number of common straight crossings on three unit intervals. Let us explain how this is done.

We first partition the interval $(0, 9)$ in disjoint intervals of length 3. This gives us the three intervals: $(0, 3)$, $(3, 6)$ and $(6, 9)$. Then we determine how many of these intervals are crossed in a straight way by $R \circ S$ when they get crossed for the first time during $a$ and and when they get crossed for the first time during $b$. In Figure 3, we see that the first crossing during $a$ of $(0, 3)$ by $R \circ S$ is straight. However, the first crossing during $b$ of $(0, 3)$ by $R \circ S$ is not. Thus, for the interval $(0, 3)$ we do not have a common first straight crossing. Next comes the interval $(3, 6)$. There, the first crossing by $R \circ S$ of $(3, 6)$ during $a$ is not straight. [That first crossing is equal to $(5, 10)$.] On the other hand, the first crossing by $R \circ S$ of $(3, 6)$ during $b$ is straight. [It is the first crossing $(t_b + 5, t_b + 8)$.] Again with the interval $(3, 6)$ we do not observe a common first straight crossing between $a$ and $b$. Eventually the first crossing by $R \circ S$ of $(6,9)$ during $a$ is straight, while the first crossing by $R \circ S$ of $(6,9)$ during $b$ is not. So, in total we have zero common straight first crossings between $a$ and $b$. When we observe few common first straight crossings between two crossings $a$ and $b$ by $S$, we decide that the crossings $a$ and $b$ took place on different places. In the example of Figure 3, the person who only observes $R \circ S$ would thus decide that the crossings $a$ and $b$ by $S$ took place on different places. In the case of Figure 4, the first crossings by $R \circ S$ of $(0,3)$ during $a$ and during $b$ are both straight. So for $(0,3)$, we have a common first straight crossing. In Figure 4 again, the first crossings by $R \circ S$ of $(3,6)$ during $a$ and during $b$ are both not straight. The first crossing by $R \circ S$ of $(6,9)$ during $a$ is straight while during $b$ it is not. Again for $(6,9)$ we do not have a common straight crossing. Thus in the case of Figure 4, the total number of "straight common first crossings" equals 1.

*General case:* Let $(t_{1a}, t_{2a})$ and $(t_{1b}, t_{2b})$ be two crossings by $R \circ S$ of the interval $(0, 3n)$. For $0 \leq m < n$, let $w_a(m)$ be equal to 1 if the first crossing by $R \circ S$ of the interval $(3m, 3m+3)$ during $(t_{1a}, t_{2a})$ is straight, and be equal to 0 otherwise. Let $w_a$ denote the binary word $w_a(0), w_a(1), w_a(2), \ldots, w_a(n-1)$. In the same manner, define the binary word $w_b$ for the crossing $(t_{1b}, t_{2b})$. The number of common straight crossings between $a$ and $b$ is defined to be the scalar product

$$w_a \times w_b := \sum_{m=0}^{n-1} w_a(m) \cdot w_b(m).$$

We use $w_a \times w_b$ as test statistic. What is its distribution under $H_0$ and under $H_1$?

Example. To have a first common straight crossing in the $H_0$-case we need three crossings to be straight while in the $H_1$-case we need four. In order to understand why this is true, look at Figure 4 first: we have there for $m = 0$ a first common straight crossing. This means that when $R \circ S$ crosses during $a$ and during $b$ for the first time $(0,3)$, we observe in both cases a straight crossing. That we have a common first straight crossing follows from the fact that the first crossing by $R$ of $(0,3)$ during $(0,13)$ is straight and the first crossings during $a$ and during $b$ of $(0,3)$ are both straight as well. In Figure 3, we have that $w_a(0) = 1$ and $w_b(0) = 0$. For $w_a(0) \cdot w_b(0)$ to be equal to 1 in Figure 1, there is only one thing missing: The first crossing $(32, 37)$ by $R$ of the interval $(0,3)$ should be straight.

*General case*: Let $m \in \mathbb{N}$ be such that $m < n$. Let $(k_{1a}, k_{2a}) = (S(t_{1a}), S(t_{2a}))$ and $(k_{1b}, k_{2b}) = (S(t_{1b}), S(t_{2b}))$. Let $(k_{1am}, k_{2am})$ designate the first crossing by $R$ of $(3m, 3m+3)$ during $(k_{1a}, k_{2a})$. Let $(k_{1bm}, k_{2bm})$ designate the first crossing by $R$ of $(3m, 3m+3)$ during $(k_{1b}, k_{2b})$. In the case of Hypothesis $H_0$ we have $(k_{1a}, k_{2a}) = (k_{1b}, k_{2b})$ and $(k_{1am}, k_{2am}) = (k_{1bm}, k_{2bm})$. We get:

*Under* $H_0$: $w_a(m) \cdot w_b(m) = 1$ iff the following *three crossings* are straight:

1. The crossing $(k_{1am}, k_{2am})$ by $R$ of the interval $(3m, 3m+3)$.
2. The first crossing by $S$ during $(t_{1a}, t_{2a})$ of the interval $(k_{1am}, k_{2am})$.
3. The first crossing by $S$ during $(t_{1b}, t_{2b})$ of the interval $(k_{1am}, k_{2am})$.

*Under* $H_1$: $w_a(m) \cdot w_b(m) = 1$ iff the following *four crossings* are straight:

1. The cr ossing $(k_{1am}, k_{2am})$ by $R$ of the interval $(3m, 3m+3)$.
2. The crossing $(k_{1b}m, k_{2bm})$ by $R$ of the interval $(3m, 3m+3)$.
3. The first crossing by $S$ during $(t_{1a}, t_{2a})$ of the interval $(k_{1am}, k_{2am})$.
4. The first crossing by $S$ during $(t_{1b}, t_{2b})$ of the interval $(k_{1bm}, k_{2bm})$.

$R$ and $S$ are independent simple random walks. For the simple random walk a crossing of an interval of length 3 is straight with probability $\frac{3}{4}$, as is shown below in Fact e.5. Under $H_0$, there are three such crossings involved, while under $H_1$ there are four. This is why $P(w_a(m) \cdot w_b(m) = 1) = (\frac{3}{4})^3$ in the case $H_0$ and $P(w_a(m) \cdot w_b(m) = 1) = (\frac{3}{4})^4$ in the case $H_1$. By the Markov property, the variables $w_a(m) \cdot w_b(m)$ for different $m$'s are independent. This gives:

The distribution of the test statistic $w_a \times w_b$ is equal to:
*Under* $H_0$: Binomial with parameter $n$ and $(\frac{3}{4})^3$.
*Under* $H_1$: Binomial with parameter $n$ and $(\frac{3}{4})^4$.
Let $c := \frac{1}{2}((\frac{3}{4})^3 + (\frac{3}{4})^4)$.
*Localization test* with parameter $n$:

(a) When $w_a \times w_b > c \cdot n$, we accept $H_0$.
(b) When $w_a \times w_b \leq c \cdot n$, we accept $H_1$.

The above statement about the distribution of the test statistic holds only if we select the pair of crossings $((t_{1a}, t_{2a}), (t_{1b}, t_{2b}))$ in an appropriate manner. For example, if we would choose $(t_{1a}, t_{2a})$ to be the first crossing by $R \circ S$ of $(0, 3n)$ such that $w_a(m) = 1$ for all $m < n$ and $(t_{1b}, t_{2b})$ to be the first crossing by $R$ of $(0, 3n)$ such that $w_a(m) = 1$ for all $m < n$, then obviously the above statement about the distributions would not hold. In Lemma 8, the statement is made rigorous. For this we need to numerate the crossings by $R \circ S$ of $(0, 3n)$, in an appropriate manner. By Lemma 5 we know that any crossing by $R \circ S$ of $(0, 3n)$ can be viewed as a crossing by $S$ of a crossing by $R$ of $(0, 3n)$. A crossing by $R \circ S$ of $(0, 3n)$ can thus be described in a unique manner as the $i$th crossing by $S$ of the $z$th crossing by $R$ of $(0, 3n)$. We index the crossings by $R$ of $(0, 3n)$ by the set $\mathbb{Z}^* := \mathbb{Z} - \{0\}$. We call the $z$th crossing by $R$ of $(0, 3n)$:

If $z > 0$, the $z$th crossing by $R(k), k \geq 0$ of $(0, 3n)$.

If $z < 0$, the $|z|$th crossing by $R(k), k \leq 0$ of $(0, 3n)$, where we count in reverse order starting at zero.

Thus, we index the crossings by $R \circ S$ of $(0, 3n)$ by the set $\mathbb{N}^* \times \mathbb{Z}^*$. For $(i, z) \in \mathbb{N}^* \times \mathbb{Z}^*$, the $(i, z)$th crossing by $R \circ S$ of $(0, 3n)$ is the crossing which corresponds to the $i$th crossing by $S$ of the $z$th crossings by $R$ of $(0, 3n)$. Picking $(t_{1a}, t_{2a})$ and $(t_{1b}, t_{2b})$ by choosing nonrandomly two elements in the index set $\mathbb{N}^* \times \mathbb{Z}^*$ makes the statement about the distribution of the test statistic rigorous. This is the content of the next lemma.

LEMMA 8. *Let $z_a, z_b \in \mathbb{Z}^*$ and let $i_a, i_b \in \mathbb{N}^*$ be nonrandom numbers. Let $(t_{1a}, t_{2a})$ and $(t_{1b}, t_{2b})$ be the two crossings by $R \circ S$ of $(0, 3n)$ for which $(t_{1a}, t_{2a})$ is the $i_a$th crossing by $S$ of the $z_a$th crossing by $R$ of $(0, 3n)$ and $(t_{1b}, t_{2b})$ is the $i_b$th crossing by $S$ of the $z_b$th crossing by $R$ of $(0, 3n)$. Then*:

*$H_0$-case [i.e., case where $z_a = z_b$ and $(S(t_{1a}), S(t_{2a})) = (S(t_{1b}), S(t_{2b}))$]: $w_a \times w_b$ has binomial distribution with parameter $n$ and $(\frac{3}{4})^3$.*

*$H_1$-case [i.e., case where $z_a \neq z_b$ and $(S(t_{1a}), S(t_{2a})) \neq (S(t_{1b}), S(t_{2b}))$]: $w_a \times w_b$ has binomial distribution with parameter $n$ and $(\frac{3}{4})^4$.*

Note that the index in $\mathbb{N}^* \times \mathbb{Z}^*$ of a crossing by $R \circ S$ of $(0, 3n)$ is not observable, (although the crossings by $R \circ S$ of $(0, 3n)$ are themselves observable). However, by large deviation for the binomial distribution, Lemma 8 guarantees that the probability of an error by our localization test is exponentially small in $n$, when the crossings compared correspond to two nonrandom indexes in $\mathbb{N}^* \times \mathbb{Z}^*$. We cannot pick crossings by their index in $\mathbb{N}^* \times \mathbb{Z}^*$ for our reconstruction algorithm, since these are not observable. Hence, the crossings we select in an observable manner have slightly different distributions from the distributions mentioned in Lemma 8. But picking the crossings in a sensible, observable manner modifies the probability of an error only slightly, so that it remains small. Next, we need to mention a few facts which are useful for the proof of Lemma 8.

FACT a. Let $M(k)_{k \in \mathbb{N}}$ be a Markov chain with state space $\mathcal{M}$. Let $a_0, a_1, a_2, \ldots$ be a sequence of (nonrandom) elements of $\mathcal{M}$. Let $\eta_{(i+1)}$ denote the first passage time of $M(k)_{k \in \mathbb{N}}$ at $a_{(i+1)}$ after $\eta_i$. Recursively: $\eta_0 := \min\{k \geq 0 | M(k) = a_0\}$. Then, $\eta_{i+1} := \min\{k \geq \eta_i | M(k) = a_{i+1}\}$. Let $Z_i$ be the path of $M$ between $\eta_i$ and $\eta_i + 1$:

$$Z_i := (M(\eta_i), M(\eta_i + 1), M(\eta_i + 2), \ldots, M(\eta_{i+1})).$$

Then, the $Z_i$'s are independent of each other.

FACT b. Let $X$ and $Z$ be two random variables living on the same space and independent of each other. Let $A$ be an event that depends only on $X$, that is, $A \in \sigma(X)$. Then conditional on $A$, $X$ and $Z$ are still independent of each other. Furthermore, conditional on $A$, $Z$ has the same marginal distribution. Thus:

$$\mathcal{L}(X, Z|A) = \mathcal{L}(X|A) \otimes \mathcal{L}(Z).$$

Fact c. Let $X_0, X_1, \ldots, X_n$ be a collection of random variables that are independent of each other. Let $A_0, A_1, \ldots, A_n$ be a collection of events such that for each $0 \leq i \leq n$, $A_i \in \sigma(X_i)$. Let $A := \bigcap_{i=0}^n A_i$. Then conditionally on $A$, the $X_i$'s are still independent of each other:

$$\mathcal{L}(X_0, X_1, \ldots, X_n|A) = \prod_{i=0}^{n} \mathcal{L}(X_i|A_i).$$

Fact d. Let $X_0, X_1, \ldots, X_n$ be a collection of random variables that are independent of each other. Let $Y_0, Y_1, \ldots, Y_n$ be a collection of random variables satisfying: conditionally on $\sigma(X_m|0 \leq m \leq n)$, the $Y_m$'s are independent of each other and their distribution depends only on their respective $X_m$'s:

$$\mathcal{L}(Y_m|X_0, X_1, \ldots, X_n) = \mathcal{L}(Y_m|X_m).$$

Let $Z_m := (X_m, Y_m)$. Then, the $Z_m$'s are independent of each other.

Fact e. Let $\kappa_0$ designate the first recurrence time of $S$ at 0, that is, $\kappa_0 := \min\{t > 0|S(t) = 0\}$. For $l > 0$, let $\kappa_l$ designate the first passage time of $S$ at $l$, that is, $\kappa_l := \min\{t|S(t) = l\}$. Let $E_{\text{cross } l}$ designate the event $\{\kappa_l < \kappa_0\}$. Let $(j_{1i}, j_{2i})$ be an increasing collection of intervals indexed by $i \in \mathbb{N}$ such that the following holds: $j_{1i} < j_{2i} \leq j_{1(i+1)}$. Assume furthermore that $j_{10} = 0$. Let $(s_{1i}, s_{2i})$ denote the first crossing by $S$ of $(j_{1i}, j_{2i})$. For natural numbers $s < t$, let $S(s,t) := (S(s), S(s+1), \ldots, S(t))$. Let $S^*(s,t)$ designate the recentered $S(s,t)$. Hence,

$$S^*(s,t) := (0, S(s+1) - S(s), \ldots, S(t) - S(s)).$$

Recall that $\Delta(s) := S(s+1) - S(s)$. Define $\Delta(s,t) := (\Delta(s), \Delta(s+1), \ldots, \Delta(t-1))$. With these definitions the following things hold:

e.1. The $S(s_{1i}, s_{2i})$'s for various $i$'s are independent of each other. Similarly, the $\Delta(s_{1i}, s_{2i})$'s are independent of each other.

Proof. Take the sequence $j_{20}, j_{21}, j_{22}, \ldots$ for the sequence $a_0, a_1, a_2, \ldots$ of Fact a. The stopping times of Fact a are then equal to $\eta_i := s_{2i}$. The crossing $(s_{1i}, s_{2i})$ happens between time $\eta_{(i-1)}$ and time $\eta_i$. By Fact a, the pieces of path of $S$ during the time intervals $[\eta_{(i-1)}, \eta_i]$ are independent of each other. Since the crossings $(s_{1i}, s_{2i})$ for different $i$'s happen during different independent time intervals, they are also independent. □

e.2. The distribution of $S^*(s_{1i}, s_{2i})$ depends only on the length $d_i := j_{2i} - j_{1i}$. The distribution of $S(s_{1i}, s_{2i})$ is equal to the distribution of the path of the random walk starting at the point $j_{1i}$ until it reaches $j_{2i}$, conditioned that it first meets $j_{2i}$ before meeting $j_{1i}$. In other words, it is conditioned on that the random walk $S$ makes a crossing of $(j_{1i}, j_{2i})$. The random walk starting at $j_{1i}$ is defined as $\{S_i(t) := S(t) + j_{1i}\}_{t\in\mathbb{N}}$. With this notation, the distribution of $S(s_{1i}, s_{2i})$ equals

$$\mathcal{L}(S_i(0, \kappa_{d_i})|E_{\text{cross } d_i})$$

or equivalently,

$$\mathcal{L}(S(t,\nu)|S(t) = j_{1i} \text{ and after time } t,\ S \text{ visits } j_{2i} \\ \text{before it returns for the first time to } j_{1i}),$$

where $t$ designates any nonrandom time, and $\nu$ designates the first visit after $t$ to $j_{2i}$.

e.3. The distribution of $\Delta(s_{1i}, s_{2i})$ depends only on the length $d_i$. It is equal to the distribution of $(\Delta(0), \Delta(1), \dots, \Delta(\kappa_{d_i}))$ conditional on the event that the random walk first meets $d_i$ before meeting 0. Thus,

$$\mathcal{L}(\Delta(s_{1i}, s_{2i})) = \mathcal{L}(\Delta(0), \Delta(1), \dots, \Delta(\kappa_{d_i})|E_{\text{cross } d_i}).$$

e.4. The joint distribution of the path of $S$ during the crossings $(s_{1i}, s_{2i})$ is not changed if we condition on the event that the crossings $(s_{1i}, s_{2i})$ have to occur during a crossing. More precisely, we are considering the joint distribution of the $(s_{1i}, s_{2i})$'s for $0 \le i \le n$. We condition under the event that we have a crossing by $S$ of $(0, j_{2n})$ starting at zero. After conditioning we get the same distribution as before:

$$\mathcal{L}(S(s_{10}, s_{20}), S(s_{11}, s_{21}), \dots, S(s_{1n}, s_{2n})) \\ = \mathcal{L}(S(s_{10}, s_{20}), S(s_{11}, s_{21}), \dots, S(s_{1n}, s_{2n})|E_{\text{cross } j_{2n}}).$$

Proof. Let $E^2_{\text{cross}}(i)$ be the event that $S$ does not visit 0 during $(s_{1i}, s_{2i})$. $E^2_{\text{cross}}(i) := \{S(t) \neq 0,\ \forall t \in (s_{1i}, s_{2i}]\}$. In a similar manner define $E^1_{\text{cross}}(i) := \{S(t) \neq 0,\ \forall t \in (s_{2(i-1)}, s_{1i}]\}$. We get

$$E_{\text{cross } j_{2n}} = \left(\bigcap_{i=0}^{n} E^1_{\text{cross}}(i)\right) \cap \left(\bigcap_{i=0}^{n} E^2_{\text{cross}}(i)\right).$$

The different pieces of paths from the collection:

$$\{S(s_{1i}, s_{2i})|0 \le i \le n\} \cup \{S(s_{1(i-1)}, s_{1i})|0 < i \le n\}$$

are independent of one another. Thus, we are exactly in the situation of Fact c. Applying Fact c to $\{S(s_{1i}, s_{2i})|0 \le i \le n\}$, we find that

$$\mathcal{L}(S(s_{10}, s_{20}), S(s_{11}, s_{21}), \dots, S(s_{1n}, s_{2n})|E_{\text{cross } j_{2n}}) \tag{1}$$

equals

$$\bigotimes_{i=0}^{n} \mathcal{L}(S(s_{1i}, s_{2i}) | E^2_{\text{cross}}(i)).$$

However, since $(s_{1i}, s_{2i})$ is a crossing by $S$ of $(j_{1i}, j_{2i})$ where $0 \le j_{1i}, j_{2i}$, it follows that a.s. $S$ during $(s_{1i}, s_{2i})$ does not visit 0. Thus the event $E^2_{\text{cross}}(i)$ is the almost sure event. Hence:

$$\mathcal{L}(S(s_{1i}, s_{2i}) | E^2_{\text{cross}}(i)) = \mathcal{L}(S(s_{1i}, s_{2i})).$$

This proves that the distribution 1 equals $\bigotimes_{i=0}^{n} \mathcal{L}(S(s_{1i}, s_{2i}))$. The last expression, by e.1, is, however, the joint distribution of the "unconditional" $S(s_{1i}, s_{2i})$'s. $\square$

e.5. The probability that a crossing by $S$ of an interval of length 3 is straight equals $\frac{3}{4}$. Thus, if $d_i = 3$, we have

$$P(s_{2i} - s_{1i} = 3) = \tfrac{3}{4}.$$

Proof. We need to calculate the probability $P(\kappa_3 = 3 | E_{\text{cross } 3})$. $E_{\text{cross } 3}$ is the event that before coming back to zero, the random walk $S$ first visits 3. It can do it in exactly $3, 5, 7, \ldots$ steps. For each given number of steps there is precisely one path. The reason is that when the random walk is in the interval $[0, 3]$, in order to not reach the border, there is always only one possible step. Any path of length $2k+1$ has probability $(\frac{1}{2})^{2k+1}$. The path of length 3 is the straight path. We find:

$$P(\kappa_3 = 3 | E_{\text{cross } 3}) = \frac{P(\kappa_3 = 3)}{P(E_{\text{cross } 3})} = \frac{(1/2)^3}{\sum_{k=1}^{\infty} (1/2)^{2k+1}} = \frac{3}{4}. \qquad \square$$

Note that Fact e holds for any simple random walk.

Fact f. Let $x_1 < x_2 \le y_1 < y_2$. Let $(t_{1xi}, t_{2xi})$ designate the $i$th crossing by $S$ of $(x_1, x_2)$. Let $(t_{1yi}, t_{2yi})$ designate the $i$th crossing by $S$ of $(y_1, y_2)$. Then, $(S(t_{1xi}, t_{2xi}))_{i \ge 0}$ is independent of $(S(t_{1yi}, t_{2yi}))_{i \ge 0}$.

Proof. Let $\iota_j$ designate the $j$th visit by $S$ to the point $x_2$. This defines a renewal process and a regenerative process. Since the random walk $S$ cannot jump, during each renewal period, it can either spend the whole time in $]\infty, x_2[$ or in $]x_2, \infty[$. During the same renewal period, $S$ cannot visit both $]\infty, x_2[$ and $]x_2, \infty[$. This implies that a crossing by $S$ of $(x_1, x_2)$ and a crossing by $S$ of $(y_1, y_2)$ can never occur during the same renewal period. The renewal periods are independent of each other; that is, the pieces of path $S(\iota_j, \iota_{j+1})$ are independent for various $j$'s. Since the crossings by $S$ of $(x_1, x_2)$ and the crossings by $S$ of $(y_1, y_2)$ occur during different independent renewal times, it follows that $(S(t_{1xi}, t_{2xi}))_{i\geq 0}$ is independent of $(S(t_{1yi}, t_{2yi}))_{i\geq 0}$. □

FACT g. Let $x_1 < x_2$ be integer numbers. Let $(t_{1xi}, t_{2xi})$ designate the $i$th crossing by $S$ of $(x_1, x_2)$. Then the pieces of path $S(t_{1xi}, t_{2xi})$ are independent of each other for various $i$'s.

PROOF. Assume without loss of generality that $0 < x_1 < x_2$. Let the sequence $a_0, a_1, a_2, \ldots$ be equal to the alternating sequence $x_1, x_2, x_1, x_2, x_1, \ldots$. Define as in Fact a the stopping times $\eta_j$. In other words, $\eta_0$ designates the first visit by $S$ to $a_0$ and $\eta_{(j+1)}$ designates the first visit by $S$ after time $\eta_j$ to the point $a_{(j+1)}$. The pieces of path in between stopping times are by Fact a independent of each other. In other words, the $S(\eta_j, \eta_{(j+1)})$'s for different $j$'s are independent. However, in each time interval $[\eta_j, \eta_{(j+1)}]$ there can be at most one crossing $(t_{1xi}, t_{2xi})$. It follows that the $S(t_{1xi}, t_{2xi})$ are independent of each other. □

NOTATION. Let $0 \leq m < n$. Let $(k_{1z_a}, k_{2z_a})$, respectively, $(k_{1z_b}, k_{2z_b})$, designate the $z_a$th, respectively, $z_b$th, crossing by $R$ of $(0, 3n)$. Let $(k_{1am}, k_{2am})$, respectively, $(k_{1bm}, k_{2bm})$, designate the first crossing by $R$ during $(k_{1z_a}, k_{2z_a})$, respectively, $(k_{1z_b}, k_{2z_b})$, of $(3m, 3m+3)$. Let $w_a^R(m)$, respectively, $w_b^R(m)$, designate the variable which is equal to 1 iff $(k_{1am}, k_{2am})$, respectively, $(k_{1bm}, k_{2bm})$, is a straight crossing. Let $(t_{1am}, t_{2am})$, respectively, $(t_{1bm}, t_{2bm})$, designate the first crossing by $S$ during $(t_{1a}, t_{2a})$, respectively, $(t_{1b}, t_{2b})$, of $(k_{1am}, k_{2am})$, respectively, $(k_{1bm}, k_{2bm})$. [Here $(t_{1a}, t_{2a})$, resp. $(t_{1b}, t_{2b})$, are defined as in Lemma 8.] Let $w_a^S(m)$, respectively, $w_b^S(m)$, designate the Bernoulli variable which is equal to 1 iff $(t_{1am}, t_{2am})$, respectively, $(t_{1bm}, t_{2bm})$, is a straight crossing. With this notation and by Lemmas 5, 6 and 7, we get $w_a^S(m) \cdot w_a^R(m) = w_a(m)$ and $w_b^S(m) \cdot w_b^R(m) = w_b(m)$. Hence, the test statistic $w_a \times w_b$ is equal to

$$\sum_{m=0}^{n-1} w_a^S(m) w_a^R(m) w_b^S(m) w_b^R(m).$$

Note that the products $w_a^S(m) w_a^R(m) w_b^S(m) w_b^R(m)$ are Bernoulli random variables. Thus to prove Lemma 8, we only need to prove that these products $w_a^S(m) w_a^R(m) w_b^S(m) w_b^R(m)$ for $m = 0, \ldots, n-1$ are i.i.d. random variables such that:

*Case* $H_0$:

$$P(w_a^S(m) w_a^R(m) w_b^S(m) w_b^R(m) = 1) = (\tfrac{3}{4})^3. \tag{2}$$

*Case* $H_1$:

$$P(w_a^S(m)w_a^R(m)w_b^S(m)w_b^R(m)=1)=(\tfrac{3}{4})^4. \tag{3}$$

PROOF OF LEMMA 8. We need to distinguish two cases:

*Case* $H_0$: In this case $z_a=z_b$ and $w_a^R(m)=w_b^R(m)$ for all $0\le m<n$. Thus,

$$w_a^S(m)w_a^R(m)w_b^S(m)w_b^R(m)=w_a^S(m)w_a^R(m)w_b^S(m).$$

It follows:

$$P(w_a(m)w_b(m)=1)=P((w_a^S(m)w_b^S(m))=1,w_a^R(m)=1).$$

The right-hand side of the last equality can be written as

$$P(w_a^S(m)w_b^S(m)=1|w_a^R(m)=1)P(w_a^R(m)=1). \tag{4}$$

We have that

$$\begin{aligned} &P(w_a^S(m)w_b^S(m)=1|w_a^R(m)=1)\\ &\qquad=E[P(w_a^S(m)w_b^S(m)=1|R(k),k\in\mathbb{Z})|w_a^R(m)=1]. \end{aligned} \tag{5}$$

The crossings $(t_{1a},t_{2a})$ and $(t_{1b},t_{2b})$ are crossings by $S$ of the random interval $(k_{1z_a},k_{2z_a})$. So Fact g does not directly apply. However, by conditioning on $\sigma(R(k),k\in\mathbb{Z})$ the interval $(k_{1z_a},k_{2z_a})$ is no longer random and we can apply Fact g: Conditioned on $\sigma(R(k),k\in\mathbb{Z})$, $S(t_{1a},t_{2a})$ and $S(t_{1b},t_{2b})$ are independent of each other. Conditional on $\sigma(R(k),k\in\mathbb{Z})$, $w_a^S(m)$ only depends on $S(t_{1a},t_{2a})$, while $w_b^S(m)$ only depends on $S(t_{1b},t_{2b})$. Hence when we condition on $R$, $w_a^S(m)$ and $w_b^S(m)$ become independent. We get

$$\begin{aligned} &P(w_a^S(m)w_b^S(m)=1|R(k),k\in\mathbb{Z})\\ &\qquad=P(w_a^S(m)=1|R(k),k\in\mathbb{Z})\cdot P(w_b^S(m)=1|R(k),k\in\mathbb{Z}). \end{aligned}$$

When $w_a^R(m)=1$, then the crossing $(k_{1am},k_{2am})$ has length 3, that is, $|k_{1am}-k_{2am}|=3$. Thus, by Facts e.4 and e.5 we find that $P(w_a^S(m)=1|w_a^R(m)=1)=\frac{3}{4}$ and $P(w_b^S(m)=1|w_a^R(m)=1)=\frac{3}{4}$. So, when $w_a^R(m)=1$ holds, we find that

$$P(w_a^S(m)w_b^S(m)=1|R(k),k\in\mathbb{Z})=(\tfrac{3}{4})^2.$$

This implies that the right-hand side of (5) is equal to $E[(\frac{3}{4})^2|w_a^R(m)=1]=(\frac{3}{4})^2$. Plugging this into (4) finishes establishing (2). Next we need to demonstrate the independence of the products $w_a^S(m)w_b^S(m)w_a^R(m)$ for $0\le m<n$ in the case $H_0$. Conditional on $\sigma(R(k),k\in\mathbb{Z})$ all of the following holds.

According to Fact g, $S(t_{1a}, t_{2a})$ is independent of $S(t_{1b}, t_{2b})$. But the $w_a^S(m)$'s for various $m$'s depend only on $S(t_{1a}, t_{2a})$ and the $w_b^S(m)$'s for various $m$'s depend only on $S(t_{1b}, t_{2b})$. Thus, $(w_a^S(m))_{0\leq m<n}$ is independent of $(w_b^S(m))_{0\leq m<n}$. Furthermore, by Fact e.1, the $w_a^S(m)$'s, respectively, the $w_b^S(m)$'s, for various $m$'s are independent of each other. This leads to that the products $w_a^S(m)w_b^S(m)$ are independent of each other. [All the last arguments were meant to hold conditionally on $\sigma(R(k), k \in \mathbb{Z})$.]

By Fact e.1, the $R(k_{1am}, k_{2am})$'s are independent among each other for various $m$'s. This puts as in the case of Fact d: Take for this $R(k_{1am}, k_{2am})$ to be $X_m$ and $Y_m$ to be $w_a^S(m)w_b^S(m)$. Conditional on $(R(k_{1am}, k_{2am}))_{0\leq m<n}$ the $w_a^S(m)w_b^S(m)$'s are independent of each other and the conditional distribution of $w_a^S(m)w_b^S(m)$ depends only on $R(k_{1am}, k_{2am})$. Fact d tells that in this case the random pairs $(w_a^S(m)w_b^S(m), R(k_{1am}, k_{2am}))$ for $0 \leq m < n$ must be independent. It follows that the products $w_a^S(m)w_b^S(m)w_a^R(m)$ are also independent of each other.

*Case $H_1$*: In this case the crossing $(k_{1z_a}, k_{2z_a})$ is different from the crossing $(k_{1z_b}, k_{2z_b})$. Fact g implies that $R(k_{1z_a}, k_{2z_a})$ is independent of $R(k_{1z_b}, k_{2z_b})$. This implies that $(R(k_{1am}, k_{2am}))_{0\leq m<n}$ is independent of $(R(k_{1bm}, k_{2bm}))_{0\leq m<n}$. Conditioned on $\sigma(R(k), k \in \mathbb{Z})$, the crossings $(t_{1a}, t_{2a})$ and $(t_{1b}, t_{2b})$ by $S$ are crossing of nonrandom intervals. Hence, conditional on $\sigma(R(k), k \in \mathbb{Z})$ and by Fact f, $S(t_{1a}, t_{2a})$ and $S(t_{1b}, t_{2b})$ are independent of one another. Fact e.2 implies that conditional on $\sigma(R(k), k \in \mathbb{Z})$, the distribution of $(S(t_{1am}, t_{2am}))_{0\leq m<n}$, respectively. $(S(t_{1bm}, t_{2bm}))_{0\leq m<n}$, depends only on $(R(k_{1am}, k_{2am}))_{0\leq m<n}$, respectively, $(R(k_{1bm}, k_{2bm}))_{0\leq m<n}$. Thus, Fact d applies, and we get that $((S(t_{1am}, t_{2am}), R(k_{1am}, k_{2am})))_{0\leq m<n}$ is independent of $((S(t_{1bm}, t_{2bm}), R(k_{1bm}, k_{2bm})))_{0\leq m<n}$. Note that $w_a^S(m)w_a^R(m)$, respectively, $w_b^S(m)w_b^R(m)$, is $\sigma((S(t_{1am}, t_{2am}), R(k_{1am}, k_{2am})))$, respectively, $\sigma((S(t_{1am}, t_{2am}), R(k_{1am}, k_{2am})))$, measurable. Thus, $(w_a^S(m)w_a^R(m))_{0\leq m<n}$ is independent of $(w_b^S(m)w_b^R(m))_{0\leq m<n}$.

Conditionally on $(R(k_{1bm}, k_{2bm}))_{0\leq m<n}$, the crossings by $S(t_{1am}, t_{2am})$ for $0 \leq m < n$ are crossings of nonrandom intervals. Hence, Fact f applies so that conditionally on $(R(k_{1bm}, k_{2bm}))_{0\leq m<n}$ the pieces of paths $S(t_{1am}, t_{2am})$ are independent of each other for various $m$'s. By Facts e.2 and e.4, conditionally on $(R(k_{1bm}, k_{2bm}))_{0\leq m<n}$, the distribution of $S(t_{1am}, t_{2am})$ depends only on $(R(k_{1am}, k_{2am}))$. However, by Fact e.1, the pieces of paths $(R(k_{1am}, k_{2am}))$ are independent of each other for various $m$'s. Thus, we can apply Fact d, and get that the pairs $(R(k_{1am}, k_{2am}), S(t_{1am}, t_{2am}))$ for $0 \leq m < n$ are independent of each other. Since $w_a^S(m)w_a^R(m)$ is $\sigma(R(k_{1am}, k_{2am}), S(t_{1am}, t_{2am}))$-measurable, it follows that the products $w_a^S(m)w_a^R(m)$ for $0 \leq m < n$ are independent of each other. In a similar way, one can show that the products $w_b^S(m)w_b^R(m)$ for $0 \leq m < n$ are independent of each other. It follows that the products $w_a^S(m)w_a^R(m)w_b^S(m)w_b^R(m)$ for various $m$'s are i.i.d. By independence of $a$ and $b$, we have that

$$P(w_a^S(m)w_a^R(m)w_b^S(m)w_b^R(m)=1)$$
$$=P(w_a^S(m)w_a^R(m)=1)P(w_b^S(m)w_b^R(m)=1).$$

The right-hand side of the last equality is equal to $P(w_a^S(m)w_a^R(m)=1)^2$, because $P(w_a^S(m)w_a^R(m)=1)=P(w_b^S(m)w_b^R(m)=1)$. Furthermore,

$$P(w_a^S(m)w_a^R(m)=1)=P(w_a^S(m)=1|w_a^R(m)=1)P(w_a^R(m)=1).$$

By Fact e.5, $P(w_a^R(m)=1)=\frac{3}{4}$. When $w_a^R(m)=1$, then $|k_{1am}-k_{2am}|=3$. $|t_{1am}-t_{2am}|$ designates the first crossing by $S$ of $(k_{1am},k_{2am})$. Thus by Fact e.5, $P(w_a^S(m)=1|w_a^R(m)=1)=\frac{3}{4}$. We are done with proving (3). □

2.4. *Details of the reconstruction algorithm.* We gave already the main ideas on how to reconstruct a finite piece of scenery. In this section we describe the technical details. Let $(k_1^{n+},k_2^{n+})$ be the first crossing after time 0 by $R$ of the interval $(0,3n)$. In other words: $k_1^{n+},k_2^{n+}\geq 0$ and for all $s,t\geq 0$ such that $(s,t)$ is a crossing by $R$ of the interval $(0,3n)$ we have $k_1^{n+}\leq s$ and $k_2^{n+}\leq t$.

Let $(k_1^{n-},k_2^{n-})$ be the last crossing before time 0 by $R$ of the interval $(0,3n)$. In other words: $k_1^{n-},k_2^{n-}\leq 0$ and for all $s,t\leq 0$ such that $(s,t)$ is a crossing by $R$ of the interval $(0,3n)$ we have $k_1^{n-}\geq s$ and $k_2^{n-}\geq t$.

In the numerical example of Figure 1, we have that $(k_1^{3+},k_2^{3+})=(0,13)$. In other words $(0,13)$ is the first crossing after zero by $R$ of $(0,9)$. The part of the graph $z\mapsto R(z)$ with $z<0$ is not represented in Figure 1, so we cannot see there $(k_1^{3-},k_2^{3-})$.

The reconstruction algorithm which reconstructs a finite piece of the scenery $\xi$, reconstructs the word $\xi(k_2^{n-}),\xi(k_2^{n-}+1),\xi(k_2^{n-}+2),\dots,\xi(k_2^{n+})$ or its transpose. It achieves this by recognizing a time interval $(r,s)$ during which the nearest-neighbor walk $S$ goes in a straight way: from the point $k_2^{n-}$ to the point $k_2^{n+}$, or from the point $k_2^{n+}$ to the point $k_2^{n-}$. $(r,s)$ is thus a straight crossing by $S$ of $(k_2^{n-},k_2^{n+})$ or of $(k_2^{n+},k_2^{n-})$. During such a straight crossing $(r,s)$ the observations reveal the piece of the scenery $\xi$ which is included between $k_2^{n-}$ and $k_2^{n+}$: $\chi(r),\chi(r+1),\chi(r+2),\dots,\chi(s)$ is equal to the word $\xi(k_2^{n-}),\xi(k_2^{n-}+1),\xi(k_2^{n-}+2),\dots,\xi(k_2^{n+})$ or its transpose. The reconstruction algorithm "for a finite piece of scenery" depends on a parameter $n$. That is why we will call it the *reconstruction algorithm at level* $n$. Thus, we have a collection of algorithms indexed by $n$. Using these algorithms for increasing $n$'s will allow us to reconstruct increasing finite pieces of the scenery $\xi$ and eventually to reconstruct the whole scenery $\xi$ up to equivalence (as a limit, after infinite time). We can already mention here that the reconstruction algorithm at level $n$ does not achieve this goal in 100% of the cases; rather, it has a small failure probability. However, this failure probability is finitely summable over $n$. This insures that only a finite number of the finite size reconstructions will contain errors. This finite number of errors has no influence on the final total reconstruction, since that one is taken to be a limit.

Next we need a few definitions and notations: let $z_1, z_2 \in \mathbb{Z}$ be such that $|z_1 - z_2|$ is a multiple of 3; that is, there exists $z \in \mathbb{Z}$ such that $z_2 - z_1 = 3z$. Let $(s_1, s_2)$ be a crossing by $R \circ S$ of $(z_1, z_2)$. Let, for $0 \leq m < |z|$, $w(m)$ be equal to 1 iff the first crossing by $R \circ S$ of $(z_1 + 3m(z/|z|), z_1 + (3m + 3)(z/|z|))$ during $(s_1, s_2)$ is straight and equal to zero otherwise. We write $w_{(s_1,s_2)}$ for the binary word:

$$w(0)w(1)w(2)\cdots w(|z|-1)$$

and call it the binary word associated with the crossing $(s_1, s_2)$ by $R \circ S$.

Among the two crossings by $R$, $(k_1^{n+}, k_2^{n+})$ and $(k_1^{n-}, k_2^{n-})$, let $(k_{1a}^n, k_{2a}^n)$ designate the one of the two which gets crossed first by $S$. In a similar way, let $(k_{1c}^n, k_{2c}^n)$ designate the other one. In this way, if $k_2^{n+}$ gets visited by $S$ before $k_2^{n-}$, we have that $(k_{1a}^n, k_{2a}^n)$ equals $(k_1^{n+}, k_2^{n+})$. Otherwise, $(k_{1a}^n, k_{2a}^n)$ equals $(k_1^{n-}, k_2^{n-})$.

Let $(t_{1i}^n, t_{2i}^n)$ designate the $i$th crossing by $R \circ S$ of the interval $(0, 3n)$. Let $w_i^n$ designate the binary words associated with the crossing $(t_{1i}^n, t_{2i}^n)$. Thus:

$$w_i^n := w_{(t_{1i}^n, t_{2i}^n)}.$$

For $z \neq 0$ with $z \in \mathbb{Z}$, let $(k_{1z}^n, k_{2z}^n)$ designate the $z$th crossing by $R$ of $(0, 3n)$. [By this we mean that if $z > 0$, then $(k_{1z}^n, k_{2z}^n)$ is the $z$th crossing after 0 by $R$ of $(0, 3n)$. If $z < 0$, $(k_{1z}^n, k_{2z}^n)$ designates the $|z|$-last crossing before 0 by $R$ of $(0, 3n)$.] Note that with this notation, we have that $(k_{11}^n, k_{21}^n) = (k_1^{n+}, k_2^{n+})$ and $(k_{1(-1)}^n, k_{2(-1)}^n) = (k_1^{n-}, k_2^{n-})$. Because $S$ starts at the origin, it cannot reach any $z$th crossing $(k_{1z}^n, k_{2z}^n)$, with $|z| > 1$, before it has not crossed $(k_1^{n+}, k_2^{n+})$ or $(k_1^{n-}, k_2^{n-})$. By Lemma 5, $(t_{11}^n, t_{21}^n)$ is also the first crossing by $S$ of a crossing by $R$ of $(0, 3n)$. It follows that $(t_{11}^n, t_{21}^n)$ is obligatorily a crossing by $S$ of either $(k_1^{n+}, k_2^{n+})$ or $(k_1^{n-}, k_2^{n-})$. Thus, $(t_{11}^n, t_{21}^n)$ is a crossing by $S$ of $(k_{1a}^n, k_{2a}^n)$.

The above discussion suggests a method for constructing stopping times which with high probability will stop the random walk at the point $k_{2a}^n$. Apply for this the localization test to the two crossings $(t_{11}^n, t_{21}^n)$ and $(t_{1i}^n, t_{2i}^n)$. If the test decides that $(t_{11}^n, t_{21}^n)$ and $(t_{1i}^n, t_{2i}^n)$ are crossings by $S$ of the same interval (i.e., Hypothesis $H_0$), decide that $S(t_{2i}^n) = k_{2a}^n$. Let $\tau^n(i)$ designate the $i$th stopping time obtained by trying to stop the random walk $S$ at $k_{2a}^n$. More precisely, $\tau^n(i)$ is equal to the $i$th, $t_{2j}^n$ for which

$$w_j^n \times w_1^n > c \cdot n.$$

The scalar product for binary words of the same length $\times$ is defined in the following way: let $w = w(0)w(1)w(2)\cdots w(k)$ and $v = v(0)v(1)v(2)\cdots v(k)$ be two binary words. $w \times v := \sum_{l=0}^{k} w(l)\cdot v(l)$. We define the relation $\leq$: $w \leq v$ iff for all $l$ with $0 \leq l \leq k$ we have that $w(l) \leq v(l)$. We define the transpose of the word $w$ and write $w^*$ for the word $w^* := w(k)w(k-1)w(k-2)\cdots w(1)$.

Let $(t^n_{1a}, t^n_{2a})$ denote the first crossing by $S$ of the interval $(k^n_{1a}, k^n_{2a})$. We have that $(t^n_{1a}, t^n_{2a}) = (t^n_{11}, t^n_{21})$. Let $(t^n_{1c}, t^n_{2c})$ denote the first crossing by $S$ of the interval $(k^n_{1c}, k^n_{2c})$. As mentioned, $(t^n_{1a}, t^n_{2a})$ is also the first crossing by $R \circ S$ of the interval $(0, 3n)$, and thus is observable. Let $w^n_a$ designate the binary word associated with the crossing $(t^n_{1a}, t^n_{2a})$ by $R \circ S$. Using our notation,

$$w^n_a := w_{(t^n_{1a}, t^n_{2a})}.$$

Note that $(t^n_{1c}, t^n_{2c})$ is also a crossing by $R \circ S$ of the interval $(0, 3n)$. Let $w_c$ denote the binary word associated with the crossing $(t^n_{1c}, t^n_{2c})$ by $R \circ S$. $(t^n_{1c}, t^n_{2c})$ and $w^n_c$ are not directly observable. We can only estimate them. We denote by $\hat{w}^n_c$ our estimate for $w^n_c$ and by $(\hat{t}^n_{1c}, \hat{t}^n_{2c})$ our estimate for $(t^n_{1c}, t^n_{2c})$. We will explain later how we obtain these estimates.

As already mentioned, the goal of the reconstruction algorithm at level $n$ is to reconstruct the finite piece of the scenery $\xi$:

$$\xi(k^n_{2c}), \xi(k^n_{2c} + u), \xi(k^n_{2c} + 2u), \ldots, \xi(k^n_{2a}).$$

[Here $u$ denotes the sign $u := (k^n_{2a} - k^n_{2c})/|(k^n_{2a} - k^n_{2c})$.] The reconstruction algorithm at level $n$ achieves this by constructing a straight crossing $(s, r)$ by $S$ of $(k^n_{2c}, k^n_{2a})$. When going from $k^n_{2c}$ to $k^n_{2a}$ in a straight way, the random walk $S$ first crosses the interval $(k^n_{2c}, k^n_{1c})$ in a straight way and then the interval $(k^n_{1a}, k^n_{2a})$. Crossing $(k^n_{2c}, k^n_{1c})$, respectively, $(k^n_{1a}, k^n_{2a})$, in a straight way, we get the maximum number of "straight crossings possible by $R \circ S$." Thus, when $(s, r)$ with $s < r$ is a straight crossing by $S$ of $(k^n_{2c}, k^n_{2a})$ we have that there exists $s_2 \leq s_1 \leq r_1 \leq r_2$ with $s_2 = s, r_2 = r$ such that $(s_2, s_1)$ is a straight crossing by $S$ of $(k^n_{2c}, k^n_{1c})$ and $(r_1, r_2)$ is a straight crossing by $S$ of the interval $(k^n_{1a}, k^n_{2a})$. In this case,

$$w_{(s_1, s_2)} \geq w^n_c \tag{6}$$

and

$$w_{(r_1, r_2)} \geq w^n_a. \tag{7}$$

The above discussion suggests a method for how to search for straight crossings $(s, r)$ by $S$ of the interval $(k^n_{2c}, k^n_{2a})$: try to find $(s, r)$ minimizing $r - s$ with $s < r$ under the following constraint: there exists $s_2 \leq s_1 \leq r_1 \leq r_2$ with $s_2 = s, r_2 = r$ such that:

1. $(s_1, s_2)$ is a crossing by $R \circ S$ of $(0, 3n)$ such that (6) is satisfied.
2. $(r_1, r_2)$ is a crossing by $R \circ S$ of $(0, 3n)$ such that (7) is satisfied.

2.5. *The reconstruction algorithm at level n.* Let $\bar{n} := n^{10.89}$ and $\dot{n} := n^{11}$. We are now ready to define the *reconstruction algorithm at level* $n$ in a precise way:

ALGORITHM 9. (i) Find $(s,r)$ minimizing $r-s$ with $s<r$ under the following constraint:

1. There exists $i \le e^{\bar{n}}$ such that $\tau^{\dot{n}}(i) \le s < r \le \tau^{\dot{n}}(i) + n^{220}$.
2. There exists $s_2 \le s_1 \le r_1 \le r_2$ with $s_2 = s, r_2 = r$ such that:
   (a) $(s_1, s_2)$ is a crossing by $R \circ S$ of $(0, 3n)$ such that $w_{(s_1,s_2)} \ge \hat{w}_c^n$ holds.
   (b) $(r_1, r_2)$ is a crossing by $R \circ S$ of $(0, 3n)$ such that $w_{(r_1,r_2)} \ge w_a^n$ holds.

(ii) The output of the reconstruction algorithm at level $n$ is the binary word which we can read in the observations $\chi$ during time $(s,r)$, that is,

$$\chi(s), \chi(s+1), \chi(s+2), \dots, \chi(r),$$

where $(s,r)$ designates the first ordered pair minimizing $r-s$ under the conditions 2(a) and 2(b).

REMARK 10. (i) $w_c^n$ is not directly observable. Thus, for our reconstruction algorithm we use the estimate $\hat{w}_c^n$ instead of $w_c^n$.

(ii) The reader might be wondering why the algorithm uses conditions 2(a) and 2(b) instead of the localization test. As a matter of fact, one could imagine to replace condition 2 by the following two conditions:

(a) $(s_1, s_2)$ is a crossing by $R \circ S$ of $(0, 3n)$ such that, when compared to the crossing $(\hat{t}_{1c}^n, \hat{t}_{2c}^n)$, the localization test decides that the two crossings occurred in the same place ($H_0$-case).

(b) $(r_1, r_2)$ is a crossing by $R \circ S$ of $(0, 3n)$ such that, when compared to the crossing $(t_{1a}^n, t_{2a}^n)$, the localization test decides that the two crossings occurred in the same place ($H_0$-case).

Replacing conditions 2(a) and 2(b) by the above conditions 1 and 2 does not work. The reason is the following: typically the points $k_{2a}^n$ and $k_{2c}^n$ are at distance order$(n^2)$ from each other. [To simplify calculations, we will just prove that the order is smaller than order$(n^9)$ and work with that.] To get at least one straight crossing by $S$ of an interval of length order$(n^2)$ we need order$(2^{n^2})$ trials. Thus our algorithm needs to be able to identify correctly order$(2^{n^2})$ crossings by $S$ of $(k_{2c}^n, k_{2a}^n)$ [in our proof order$(2^{n^9})$]. The localization algorithm (with parameter $n$) has a positive probability of making an error of order$(e^{-k\cdot n})$ where $k>0$ is a constant not depending on $n$. With order$(2^{n^9})$ trials we can be sure that the localization test (with parameter $n$) will make many errors, and thus cannot be used instead of conditions 2(a) and 2(b).

(iii) If we perform the localization test with parameter $\dot{n}$ instead of $n$, the probability of an error is of order$(e^{-k\dot{n}})$. This is so small that, with high probability, we can apply it order$(e^{k\bar{n}})$ times without making a single mistake. This is more than enough trials to get, with high probability, one straight crossing by $S$ of an interval of length order$(n^9)$. This is why for condition 1 in the reconstruction algorithm at level $n$, we construct the stopping times $\tau^{\dot{n}}(i)$ using the localization algorithm with parameter $\dot{n}$.

(iv) The conditions 2(a) and 2(b) can be seen as a modified version of the localization algorithm with parameter $n$. We will show that with high probability within distance $n^{220}$ of the point $k_{2a}^{\dot{n}}$ we have: only the crossing $(k_{1a}^n, k_{2a}^n)$ is such that a crossing $(r_1, r_2)$ by $S$ of it can satisfy the inequality $w_{(r_1,r_2)} \le w_a$. A similar condition also holds for $(k_{1c}^n, k_{2c}^n)$. This implies that as long as we are within distance $n^{220}$ of the point $k_{2a}^{\dot{n}}$, conditions 2(a) and 2(b) can never make a mistake at identifying crossings by $S$ of $(k_{1a}^n, k_{2a}^n)$ and of $(k_{1c}^n, k_{2c}^n)$. When $S(\tau^{\dot{n}}(i)) = k_{2a}^{\dot{n}}$, then by definition, a crossing $(s, r)$ satisfying condition 1 of the selection rule of the reconstruction algorithm at level $n$, is such that $S(s)$ and $S(r)$ are within distance $n^{220}$ of the point $k_{2a}^{\dot{n}}$. For more details about why the reconstruction algorithm at level $n$ works, see Section 4.

2.6. *Construction of* $(\hat{t}_{1c}^n, \hat{t}_{2c}^n)$ *and of* $\hat{w}_c^n$. Recall that a crossing $(s,t)$ is called positive if $s < t$ and negative otherwise. Recall also that from the two crossings $(k_{11}^n, k_{21}^n)$ and $(k_{1(-1)}^n, k_{2(-1)}^n)$ by $R$ of $(0, 3n)$, the one which gets first crossed by $S$ is called $(k_{1a}^n, k_{2a}^n)$ while the other one is called $(k_{1c}^n, k_{2c}^n)$. After having crossed from the point $k_{1a}^n$ to the point $k_{2a}^n$, $S$ first needs to cross back from the point $k_{2a}^n$ to the point $k_{1a}^n$ before being able to cross $(k_{1c}^n, k_{2c}^n)$. More precisely, after a positive crossing by $S$ of $(k_{1a}^n, k_{2a}^n)$ there first needs to be a negative crossing by $S$ of $(k_{1a}^n, k_{2a}^n)$ before there can be a crossing by $S$ of $(k_{1c}^n, k_{2c}^n)$. On the other hand, right after a negative crossing by $S$ of $(k_{1a}^n, k_{2a}^n)$ the random walk $S$ is always located between the points $k_{1a}^n$ and $k_{1c}^n$. When the random walk $S$ is located between $k_{1a}^n$ and $k_{1c}^n$, the next time it crosses an interval $(k_{1z}^n, k_{2z}^n)$ this must be the interval $(k_{1a}^n, k_{2a}^n)$ or $(k_{1c}^n, k_{2c}^n)$. This gives a way to characterize $(t_{1c}^n, t_{2c}^n)$ [recall that $(t_{1c}^n, t_{2c}^n)$ is the first crossing by $S$ of $(k_{1c}^n, k_{2c}^n)$]: $(t_{1c}^n, t_{2c}^n)$ is the first crossing by $S$ of an interval $(k_{1z}^n, k_{2z}^n)$ such that the following two conditions are satisfied:

(i) $(t_{1c}^n, t_{2c}^n)$ is not a crossing by $S$ of $(k_{1a}^n, k_{2a}^n)$.

(ii) the last crossing by $S$ of an interval $(k_{1z}^n, k_{2z}^n)$ before $(t_{1c}^n, t_{2c}^n)$ is a negative crossing by $S$ of $(k_{1a}^n, k_{2a}^n)$.

Note that Lemma 5 implies that the crossings by $S$ of an interval $(k_{1z}^n, k_{2z}^n)$ can be characterized as follows: $(s,t)$ is a crossing by $S$ of an interval $(k_{1z}^n, k_{2z}^n)$ iff $(s,t)$ is a crossing by $R \circ S$ of $(0, 3n)$. Applying the last characterization to the above conditions leads to: $(t_{1c}^n, t_{2c}^n)$ is equal to the first crossing $(t^n_{1i}, t^n_{2i)})$ by $R \circ S$ of $(0, 3n)$ with $i > 1$ such that the following two conditions hold:

(i) $(t^n_{1i}, t^n_{2i})$ is not a crossing by $S$ of $(k^n_{1a}, k^n_{2a})$.
(ii) $(t^n_{1(i-1)}, t^n_{2(i-1)})$ is a negative crossing by $S$ of $(k^n_{1a}, k^n_{2a})$.

Which crossings are crossings by $R \circ S$ of $(0, 3n)$ is observable. That means that the crossings $(t^n_{1i}, t^n_{2i})$ are known to us. On the other hand, which crossings are crossings by $S$ of $(k^n_{1a}, k^n_{2a})$ is not directly observable. However, $(t^n_{11}, t^n_{21})$ is observable and is a crossing by $S$ of $(k^n_{1a}, k^n_{2a})$. So we can estimate if $(t^n_{1i}, t^n_{2i})$ is a crossing by $S$ of $(k^n_{1a}, k^n_{2a})$ or not. For this we ask our localization test to compare the crossings $(t^n_{11}, t^n_{21})$ and $(t^n_{1i}, t^n_{2i})$. The localization test can then estimate if the crossings $(t^n_{11}, t^n_{21})$ and $(t^n_{1i}, t^n_{2i})$ of $S$ occur on the same place or not. Our estimate for $(t^n_{1c}, t^n_{2c})$ will be defined to be the first $(t^n_{1i}, t^n_{2i})$ for which the above characterizing conditions are estimated to be true:

We define $(\hat{t}^n_{1c}, \hat{t}^n_{2c})$ to be equal to the first $(t^n_{1i}, t^n_{2i})$ with $i > 1$ for which the following three conditions hold:

(i) The localization test, when comparing $(t^n_{11}, t^n_{21})$ with $(t^n_{1i}, t^n_{2i})$, rejects the $H_0$-hypothesis.
(ii) $t^n_{1(i-1)} > t^n_{2(i-1)}$.
(iii) The localization test, when comparing $(t^n_{11}, t^n_{21})$ with $(t^n_{1(i-1)}, t^n_{2(i-1)})$, accepts the $H_0$-hypothesis.

We define $\hat{w}^n_c$ to be the binary word associated with the crossing $(\hat{t}^n_{1c}, \hat{t}^n_{2c})$.

**3. Assembling the pieces.** The reconstruction algorithm at level $n$ tries to reconstruct the finite piece of the scenery $\xi$:

$$\xi^n := \xi(k^n_{1c}), \xi(k^n_{1c} + u), \xi(k^n_{1c} + 2u), \ldots, \xi(k^n_{1a}),$$

where $u := (k^n_{1a} - k^n_{1c})/|(k^n_{1a} - k^n_{1c})|$. In this section, we explain how to construct a scenery $\bar{\xi} : \mathbb{Z} \to \{0, 1\}$, equivalent to $\xi$, from the collection of finite pieces $\xi^1, \xi^2, \ldots$. The reconstruction algorithm at level $n$ gives us the binary word $\xi^n$, but does not tell us where it is located in the scenery $\xi$. This implies that we need to "assemble" the pieces $\xi^n$ in order to get $\bar{\xi}$.

Let us introduce a few definitions: let $v = v(0)v(1)v(2) \cdots v(i)$ and $w = w(0)w(1)w(2) \cdots w(j)$ be two binary words. We say that $v$ is contained in $w$ iff there exist $j_1, j_2 \in \{0, 1, 2, \ldots, j\}$ such that $v$ is equal to

$$v = w(j_1)w(j_1 + u)w(j_1 + 2u) \cdots w(j_2) \tag{8}$$

where $u := (j_2 - j_1)/|j_2 - j_1|$. We write $v \preccurlyeq w$ when $v$ is contained in $w$. We say that $v$ is uniquely contained in $w$ and write $v \preccurlyeq_1 w$, iff there exists exactly one ordered pair $(j_1, j_2)$ in $\{0, 1, 2, \ldots, j\}^2$ such that (8) is satisfied.

Note that the sequence of pieces $\xi^1,\xi^2,\dots$ is an increasing sequence, in the sense that $\xi^n \preccurlyeq \xi^{n+1}$ for all $n\in\mathbb{N}$. (The reason for this being true is that by definition, $k_2^{n-} > k_2^{(n+1)-}$ and $k_2^{n+} < k_2^{(n+1)+}$ for all $n\in\mathbb{N}$. Thus the interval with the two endpoints $k_{2c}^n, k_{2a}^n$ is contained in the interval with endpoints $k_{2c}^m, k_{2a}^m$ when $n<m$.) Imagine that not only $\xi^n \preccurlyeq \xi^{n+1}$, but even $\xi^n \preccurlyeq_1 \xi^{n+1}$ for all $n\in\mathbb{N}$. Then there would be a unique way to assemble the pieces $\xi^1,\xi^2,\xi^3,\dots$. The situation in this case is similar to that of a puzzle: for a puzzle, once we have decided on the position of one piece, there is a unique way to assemble the whole puzzle. Furthermore, when we assemble a puzzle we always get the same image up to an isometric mapping. This is exactly the situation we encounter with the pieces of scenery when $\xi^n \preccurlyeq_1 \xi^{n+1}$ for all $n\in\mathbb{N}$.

Let us illustrate this with a practical example. Let $\xi:\mathbb{Z}\to\{0,1\}$ be the scenery from which we show below a finite portion close to the origin:

| $\xi(k)$ | ... | 1 | 0 | 1 | 0 | 0 | 0 | 1 | 1 | 1 | 0 | 0 | 1 | ... |
|---|---|---|---|---|---|---|---|---|---|---|---|---|---|---|
| $k$ | ... | $-4$ | $-3$ | $-2$ | $-1$ | 0 | 1 | 2 | 3 | 4 | 5 | 6 | 7 | ... |

Assume that we would be given the three pieces (of the part of the scenery $\xi$ which is represented above): 11000, 1000111 and 0100011100. In this case the first piece lies in the scenery $\xi$ between the points 3 and $-1$. The second piece is the piece of $\xi$ which lies between $-1$ and 4. The last piece lies between the points $-3$ and 6. We see that the first piece is uniquely contained in the second which itself is uniquely contained in the third piece. To assemble the three pieces we first place the first piece anywhere in $\mathbb{Z}$. Then we place the second piece so that it covers the first piece, and so that on the first piece it coincides with the first piece. Eventually we place the third piece so that it coincides with and covers the second one. If we place the first piece starting at the origin we get

| $\bar{\xi}(k)$ | | | | | 1 | 1 | 0 | 0 | 0 | | | |
|---|---|---|---|---|---|---|---|---|---|---|---|---|
| $k$ | $-4$ | $-3$ | $-2$ | $-1$ | 0 | 1 | 2 | 3 | 4 | 5 | 6 | 7 |

Then we place the second piece so that it covers and coincides with the first piece. For this we have to turn the second piece around. We obtain

| $\bar{\xi}(k)$ | | | | 1 | 1 | 1 | 0 | 0 | 0 | 1 | | |
|---|---|---|---|---|---|---|---|---|---|---|---|---|
| $k$ | $-4$ | $-3$ | $-2$ | $-1$ | 0 | 1 | 2 | 3 | 4 | 5 | 6 | 7 |

Eventually we place the third word and get

| $\bar{\xi}(k)$ | | 0 | 0 | 1 | 1 | 1 | 0 | 0 | 0 | 1 | 0 | |
|---|---|---|---|---|---|---|---|---|---|---|---|---|
| $k$ | $-4$ | $-3$ | $-2$ | $-1$ | 0 | 1 | 2 | 3 | 4 | 5 | 6 | 7 |

If we would go on with more and more pieces, as $n$ tends to infinity we would obtain a scenery $\bar{\xi}$ which is equivalent to $\xi$.

Let $E_0^n$ denote the event that

$$E_0^n = \{\xi^n \preccurlyeq_1 \xi^{n+1}\}.$$

We will show that

$$\sum_{n=1}^{\infty} P(E_0^{nc}) < \infty,$$

where $E_0^{nc}$ denotes the complement of $E_0^n$. From the last inequality above it follows that a.s. for all but a finite number of $n$'s we have that $\xi^n \preccurlyeq_1 \xi^{n+1}$. The assemblage procedure we define below still works if $\xi^n \preccurlyeq_1 \xi^{n+1}$ holds for all but a finite number of $n$'s.

Let us mention an additional problem: each reconstruction algorithm at level $n$ has a small probability of making an error. Thus the output of the reconstruction algorithm at level $n$ is not a.s. equal to $\xi^n$ but is only an estimate of $\xi^n$. For the output of the reconstruction algorithm at level $n$, we will thus write $\hat{\xi}^n$ instead of $\xi^n$. We denote by $E^n$ the event that the algorithm at level $n$ works. That is,

$$E^n := \{\xi^n = \hat{\xi}^n\}.$$

By $E^{nc}$ we denote the complementary event of $E^n$. In the next section it is shown that

$$\sum_{n=1}^{\infty} P(E^{nc}) < \infty. \tag{9}$$

From this it follows that almost surely all but a finite number of reconstructions $\hat{\xi}^n$ are correct, that is, are such that $\xi^n = \hat{\xi}^n$. Our assembling procedure defined below is robust against this kind of problem: if only a finite number of pieces $\hat{\xi}^n$ are wrong it still works. Let us next define in a precise way our *assemblage procedure*:

ALGORITHM 11. (i) Let $l^n + 1$ designate the length of the word $\hat{\xi}^n$ and let $\hat{\xi}^n(i)$ designate the $i$th bit of the binary word $\hat{\xi}^n$. In this way,

$$\hat{\xi}^n = \hat{\xi}^n(0)\hat{\xi}^n(1)\hat{\xi}^n(2)\cdots\hat{\xi}^n(l^n).$$

(ii) Let $n_0$ designate the smallest natural (random) number such that for all $n \geq n_0$ we have that $\hat{\xi}^n \preccurlyeq_1 \hat{\xi}^{n+1}$ holds.

(iii) We construct the scenery $\bar{\xi}$ by induction on $n$ starting at $n_0$.

We first place the word $\hat{\xi}^{n_0}$ at the origin.

Once the word $\hat{\xi}^n$ is placed, we place the word $\hat{\xi}^{n+1}$ in the unique manner such that it covers and coincides with $\hat{\xi}^n$ on $\hat{\xi}^n$.

$(d_1^n, d_2^n)$ designates the position of where we placed the word $\hat{\xi}^n$. More precisely,

(a) Let $d_1^{n_0} := 0$ and let $d_2^{n_0} := l^{n_0}$. For all $k \in [0, d_2^{n_0}]$ define: $\bar{\xi}(k) := \hat{\xi}^{n_0}(k)$.

(b) Once $d_1^n, d_2^n$ are defined and $\bar{\xi}(k)$ is defined for all $k \in [d_1^n, d_2^n]$ let $d_1^{n+1}, d_2^{n+1}$ with $d_1^{n+1} \leq d_2^{n+1}$ be the unique ordered pair of integers such that $[d_1^n, d_2^n] \subset [d_1^{n+1}, d_2^{n+1}]$ and such that one of the following two cases holds:

1. For all $k \in [d_1^n, d_2^n]$ we have that
$$\bar{\xi}(k) = \hat{\xi}^{n+1}(k - d_1^{n+1}).$$
2. For all $k \in [d_1^n, d_2^n]$ we have that
$$\bar{\xi}(k) = \hat{\xi}^{n+1}(l^{n+1} - (k - d_1^{n+1})).$$

For all $k \in [d_1^{n+1}, d_2^{n+1}]$, let $\bar{\xi}(k)$ be equal to

1. When case 1 above holds,
$$\bar{\xi}(k) := \hat{\xi}^{n+1}(k - d_1^{n+1}).$$
2. When case 2 above holds,
$$\bar{\xi}(k) := \hat{\xi}^{n+1}(l^{n+1} - k - d_1^{n+1}).$$

The constructed scenery $\bar{\xi}$ is equivalent to $\xi$ as soon as for all but a finite number of $n$'s we have that $\xi^n \preccurlyeq_1 \xi^{n+1}$ and $\xi^n = \hat{\xi}^n$. This should be obvious and we leave the proof to the reader. It thus only remains to prove that almost surely for all but a finite number of $n$'s, $\xi^n \preccurlyeq_1 \xi^{n+1}$ and $\xi^n = \hat{\xi}^n$ hold.

**4. Proof that the reconstruction at level $n$ works.** In this section we prove that the reconstruction algorithm at level $n$ works with high probability; that is, we prove (9). For this we decompose $E^n$ into several elementary events. Let us start with some definition.

We say that $(s, r)$ satisfies the conditions of Algorithm 9 with $w_c^n$ instead of $\hat{w}_c^n$ iff $s < r$ and it satisfies all of the following conditions:

1. There exists $i \leq e^{\bar{n}}$ such that $\tau^{\dot{n}}(i) \leq s < r \leq \tau^{\dot{n}}(i) + n^{220}$.
2. There exists $s_2 \leq s_1 \leq r_1 \leq r_2$ with $s_2 = s, r_2 = r$ such that:
   (a) $(s_1, s_2)$ is a crossing by $R \circ S$ of $(0, 3n)$ such that $w_{(s_1,s_2)} \geq w_c^n$ holds.
   (b) $(r_1, r_2)$ is a crossing by $R \circ S$ of $(0, 3n)$ such that $w_{(r_1,r_2)} \geq w_a^n$.

Let $E_1^n$ designate the event that if Algorithm 9 is given the real $w_c^n$ instead of the estimate $\hat{w}_c^n$, it produces a straight crossing by $S$ of $(k_{2c}^n, k_{2a}^n)$:

$$E_1^n := \{\text{There exists at least one pair } (s, r), \text{ satisfying the conditions of Algorithm 9 with } w_c^n \text{ instead of } \hat{w}_c^n\}$$
$$\cap \{\text{Any pair } (s, r), \text{ minimizing } r - s \text{ under the conditions of Algorithm 9 with } w_c^n \text{ instead of } \hat{w}_c^n, \text{ is a straight crossing by } S \text{ of } (k_{2c}^n, k_{2a}^n)\}.$$

Let $E^n_{t_c}$ be the event that the construction of $(t^n_{c1}, t^n_{c2})$ works:

$$E^n_{t_c} := \{(\hat{t}^n_{c1}, \hat{t}^n_{c2}) = (t^n_{c1}, t^n_{c2})\}. \tag{10}$$

Note that when $E^n_{t_c}$ holds, then $w^n_c = \hat{w}^n_c$:

$$\begin{aligned} E^n_{\text{all correct}} := \{&\text{All } (s,r) \text{ satisfying the constraints of Algorithm 9} \\ &\text{with } w^n_c \text{ instead of } \hat{w}^n_c, \text{ are such that } S(s) = k^n_{2c}, S(r) = k^n_{2a}\}, \\ E^n_{\text{at least one}} := \{&\text{There exists } (s,r) \text{ satisfying the constraints of Algorithm 9} \\ &\text{with } w^n_c \text{ instead of } \hat{w}^n_c, \text{ such that } (s,r) \text{ is a straight crossing} \\ &\text{by } S \text{ of } (k^n_{2c}, k^n_{2a})\}. \end{aligned}$$

Let $(t^n_{1ai}, t^n_{2ai})$ be the $i$th crossing by $S$ of $(k^n_{1a}, k^n_{2a})$. Let $E^n_{\text{stopping}}$ be the event that the stopping times $\tau^n(i)$ stop the random walk at $k^n_{2a}$:

$$E^n_{\text{stopping}} := \{t^n_{2ai} = \tau^n(i),\ \forall\, i \le \exp(n^{0.99})\}.$$

Let:

$$\begin{aligned} &E^n_{\text{no other } a \text{ crossing by } R} \\ &\quad := \{\text{The only crossing } (k_1,k_2) \text{ by } R \text{ of } (0,3n) \\ &\qquad \text{with } |k_1 - k^{\dot{n}}_{2a}|, |k_2 - k^{\dot{n}}_{2a}| \le n^{220} \text{ such that } w^R_{(k_1,k_2)} \ge w^n_a \text{ is } (k^n_{1a}, k^n_{2a})\}, \\ &E^n_{\text{no other } c \text{ crossing by } R} \\ &\quad := \{\text{The only crossing } (k_1,k_2) \text{ by } R \text{ of } (0,3n) \\ &\qquad \text{with } |k_1 - k^{\dot{n}}_{2a}|, |k_2 - k^{\dot{n}}_{2a}| \le n^{220} \text{ such that } w^R_{(k_1,k_2)} \ge w^n_c \text{ is } (k^n_{1c}, k^n_{2c})\}, \\ &E^n_{\text{no other crossing by } R} \\ &\quad := E^n_{\text{no other } a \text{ crossing by } R} \cap E^n_{\text{no other } c \text{ crossing by } R}, \\ &E^n_{\text{straight}} \\ &\quad := \{\text{There exists } i \le e^{\bar{n}} \text{ and } s,r \text{ with } t^{\dot{n}}_{2ai} \le s,r \le t^{\dot{n}}_{2ai} + n^{220} \\ &\qquad \text{such that } (s,r) \text{ is a straight crossing by } S \text{ of } (k^n_{1c}, k^n_{2a})\}. \end{aligned}$$

Let $E^n_{\text{visit}}$ be the event that the random walk $S$ visits the point $k^n_{2c}$ before time $\exp(n^{0.5})$:

$$E^n_{\text{visit}} := \{t^n_{2c} < \exp(n^{0.5})\}.$$

Recall that $\dot{n} := n^{11}$. In Section 4.1 we prove the following inclusions:

$$E^n_1 \cap E^n_{t_c} \subset E^n, \tag{11}$$

$$E^n_{\text{at least one}} \cap E^n_{\text{all correct}} \subset E^n_1, \tag{12}$$

$$E^{\dot n}_{\text{stopping}} \cap E^n_{\text{no other crossing by } R} \subset E^n_{\text{all correct}}, \tag{13}$$

$$E^n_{\text{straight}} \cap E^{\dot n}_{\text{stopping}} \subset E^n_{\text{at least one}}, \tag{14}$$

$$E^n_{\text{stopping}} \cap E^n_{\text{visit}} \subset E_{t_c}. \tag{15}$$

From the inclusions (11)–(15) it follows that

$$E^n_{\text{straight}} \cap E^{\dot n}_{\text{stopping}} \cap E^n_{\text{stopping}} \cap E^n_{\text{no other crossing by } R} \cap E^n_{\text{visit}} \subset E^n,$$

which implies

$$P(E^{nc}_{\text{straight}}) + P(E^{\dot n c}_{\text{stopping}}) + P(E^{nc}_{\text{stopping}}) \\ + P(E^{nc}_{\text{no other crossing by } R}) + P(E^{nc}_{\text{visit}}) \geq P(E^{nc}).$$

(Here $E^{nc}_{\text{something}}$ designates the complement of the event $E^n_{\text{something}}$.) In Section 4.2 we prove that

$$P(E^{nc}_{\text{straight}}), \qquad P(E^{\dot n c}_{\text{stopping}}), \qquad P(E^{nc}_{\text{stopping}}),$$
$$P(E^{nc}_{\text{no other crossing by } R}) \quad \text{and} \quad P(E^{nc}_{\text{visit}})$$

are all finitely summable over $n$. Together with the last inequality, this proves that $P(E^{nc})$ is finitely summable over $n$.

### 4.1. *Combinatorics.*

Proof that $E^n_1 \cap E^n_{t_c} \subset E^n$ holds. When $E^n_{t_c}$ holds, then $w^n_c = \hat{w}^n_c$. In this case, the event $E^n_1$ amounts to the same as event $E^n$. It follows that $E^n_1 \cap E^n_{t_c} = E^n \cap E^n_{t_c}$, which implies inclusion (11). □

Proof that $E^n_{\text{at least one}} \cap E^n_{\text{all correct}} \subset E^n_1$ holds. Let $(s, r)$ be a pair minimizing $r - s$ under the constraint of Algorithm 9 with $w^n_c$ instead of $\hat{w}^n_c$. Then if $E^n_{\text{all correct}}$ holds, we have that $S(s) = k^n_{2c}, S(r) = k^n_{2a}$. If $E^n_{\text{at least one}}$ also holds, there exists a straight crossing $(s', r')$ by $S$ of $(k^n_{2c}, k^n_{2a})$ satisfying the constraint of Algorithm 9 with $w^n_c$ instead of $\hat{w}^n_c$. For a straight crossing we have $r' - s' = |k^n_{2c} - k^n_{2a}|$. Since $r - s$ is minimal under the constraint of Algorithm 9, we get $|r - s| \leq |k^n_{2c} - k^n_{2a}|$. This together with $S(s) = k^n_{2c}, S(r) = k^n_{2a}$ is only possible if $(s, r)$ is a straight crossing by $S$ of $(k^n_{2c}, k^n_{2a})$. We just proved that when $E^n_{\text{at least one}}$ and $E^n_{\text{all correct}}$ hold, all pair $(s, r)$ minimizing $r - s$ under the constraint of Algorithm 9 are straight crossings by $S$ of $(k^n_{2c}, k^n_{2a})$. In this case $E^n_1$ holds. Thus, together $E^n_{\text{at least one}}$ and $E^n_{\text{all correct}}$ imply $E^n_1$. □

Proof that $E^{\dot n}_{\text{stopping}} \cap E^n_{\text{no other crossing by } R} \subset E^n_{\text{all correct}}$ holds. Let $(s, r)$ satisfy all the constraints of Algorithm 9. Then there exists $s_2 \leq s_1 \leq r_1 \leq r_2$ with $s_2 = s, r_2 = r$, where $(r_1, r_2)$ is a crossing by $R \circ S$ of $(0, 3n)$ such that $w_{(r_1,r_2)} \geq w_a^n$ holds. By Lemma 5 we have that there exists a crossing $(k_1,k_2)$ by $R$ of $(0,3n)$ such that $(r_1,r_2)$ is a crossing by $S$ of $(k_1,k_2)$. By Lemmas 6 and 7, we have that $w^R_{(k_1,k_2)} \geq w_{(r_1,r_2)}$. Thus, $w^R_{(k_1,k_2)} \geq w_a^n$.

Additionally by the constraints of Algorithm 9 there exists $i \leq e^{\bar{n}}$ such that $\tau^{\dot{n}}(i) \leq s < r \leq \tau^{\dot{n}}(i)+n^{220}$. If additionally $E^{\dot{n}}_{\text{stopping}}$ holds, then $S(\tau^{\dot{n}}(i)) = k^{\dot{n}}_{2a}$. The random walk $S$ during a time interval of $n^{220}$ time cannot walk further than $n^{220}$. Thus, $|S(r_1) - k^{\dot{n}}_{2a}|, |S(r_2) - k^{\dot{n}}_{2a}| \leq n^{220}$. This is equivalent to saying that $|k_1 - k^{\dot{n}}_{2a}|, |k_2 - k^{\dot{n}}_{2a}| \leq n^{220}$. Hence the condition in event $E^n_{\text{no other crossing by } R}$ applies to the crossing $(k_1,k_2)$. It follows that if $E^n_{\text{no other crossing by } R}$ also holds, then $(k_1,k_2)$ equals $(k^n_{1a}, k^n_{2a})$. This implies that $S(r) = k^n_{2a}$. We have proven that when $E^{\dot{n}}_{\text{stopping}}$ and $E^n_{\text{no other crossing by } R}$ both hold, then $S(r) = k^n_{2a}$. In a similar way, one can prove that in this case $S(s) = k^n_{2c}$. (We leave that proof to the reader.) Thus, $E^{\dot{n}}_{\text{stopping}}$ and $E^n_{\text{no other crossing by } R}$ jointly imply $E^n_{\text{all correct}}$. □

Proof that $E^n_{\text{straight}} \cap E^{\dot{n}}_{\text{stopping}} \subset E_{\text{at least one}}$ holds. $E^n_{\text{straight}}$ and $E^{\dot{n}}_{\text{stopping}}$ jointly imply that there exist $i \leq e^{\bar{n}}$ and $s, r$ with $\tau^{\dot{n}}(i) \leq s, r \leq \tau^{\dot{n}}(i) + n^{220}$ such that $(s,r)$ is a straight crossing by $S$ of $(k^n_{2c}, k^n_{2a})$. Thus, $(s,r)$ already satisfies condition 1 of Algorithm 9. It remains to show that $(s,r)$ also satisfies condition 2. During the time interval $(s,r)$, $S$ crosses from the point $k^n_{2c}$ to the point $k^n_{2a}$ in a straight way. For this, $S$ first needs to cross $(k^n_{2c}, k^n_{1c})$ in a straight manner and then $(k^n_{1a}, k^n_{2a})$. Thus, there exists $s_2 \leq s_1 \leq r_1 \leq r_2$ with $s_2 = s, r_2 = r$ such that $(s_2, s_1)$ is a straight crossing by $S$ of $(k^n_{2c}, k^n_{1c})$ and $(r_1, r_2)$ is a straight crossing by $S$ of $(k^n_{1a}, k^n_{2a})$. We know by Lemma 5 that a crossing of a crossing is a crossing of the composition. Thus, $(s_1, s_2)$ and $(r_1, r_2)$ are both crossings by $R \circ S$ of $(0,3n)$. Since the crossing $(s_1,s_2)$ by $S$ is straight, we have by Lemmas 6 and 7 that $w_{(s_1,s_2)} = w^R_{(k^n_{2c},k^n_{1c})}$. By Lemmas 6 and 7 again, we have that $w^R_{(k^n_{2c},k^n_{1c})} \geq w^n_c$. Thus, $w_{(s_1,s_2)} \geq w^n_c$. In a similar way one can show that $w_{(r_1,r_2)} \geq w^n_a$. This proves that $(s,r)$ satisfies the conditions of Algorithm 9 with $w^n_c$ instead of $\hat{w}^n_c$. However, $(s,r)$ is a straight crossing by $S$ of $(k^n_{2c}, k^n_{2a})$. Thus, $E_{\text{at least one}}$ holds. We just proved that $E^n_{\text{straight}}$ and $E^{\dot{n}}_{\text{stopping}}$ together imply $E_{\text{at least one}}$. □

Proof that $E^n_{\text{stopping}} \cap E^n_{\text{visit}} \subset E_{t_c}$ holds. In Section 2.6, we saw that $(t^n_{1c}, t^n_{2c})$ can be characterized as follows: $(t^n_{1c}, t^n_{2c})$ is equal to the first crossing $(t^n_{1i}, t^n_{2i})$ by $R \circ S$ of $(0,3n)$ with $i > 1$ such that the following two conditions hold:

(i) $(t^n_{1i}, t^n_{2i})$ is not a crossing by $S$ of $(k^n_{1a}, k^n_{2a})$.
(ii) $(t^n_{1(i-1)}, t^n_{2(i-1)})$ is a negative crossing by $S$ of $(k^n_{1a}, k^n_{2a})$.

The estimate $(\hat{t}^n_{c1}, \hat{t}^n_{c2})$ is defined to be the first crossing by $R \circ S$ of $(0,3n)$ for which our localization test decides that the two conditions in the last characterization above hold. Thus, if up to time $t^n_{c2}$, the localization test gets all the crossings by $S$ of $(k^n_{1a}, k^n_{2a})$ right, then the reconstruction of $(t^n_{1c}, t^n_{2c})$ works, that is, $E_{t_c}$ holds. The event $E^n_{\text{stopping}}$ tells us that up to $t^n_{2ai}$ with $i = \exp(n^{0.99})$ the localization test makes no errors in recognizing the crossings by $S$ of $(k^n_{1a}, k^n_{2a})$. However, $\exp(n^{0.99}) \leq t^n_{2ai}$ for $i = \exp(n^{0.99})$, since each crossing lasts at least one time unit. Thus, up to time $\exp(n^{0.99})$ the localization test makes no errors in recognizing the crossings by $S$ of $(k^n_{1a}, k^n_{2a})$. However, if $E^n_{\text{visit}}$ holds, then the random walk $S$ visits the point $k^n_{2c}$ before time $\exp(n^{0.5})$. Also, $\exp(n^{0.5}) \leq \exp(n^{0.99})$. In that case, before $S$ visits the point $k^n_{2c}$, no errors occur. This proves that $E^n_{\text{stopping}}$ and $E^n_{\text{visit}}$ jointly imply $E_{t_c}$. $\square$

4.2. *Probability bounds.*

*High probability of $E^n_{\text{visit}}$.* We need a few definitions. Let $E^n_2$ be the event that the random walk $S$ visits both points $n^{10}$ and $-n^{10}$ before time $\exp(n^{0.5})$. Let

$$E^n_{k_a,c} := \{|k^n_{2a}|, |k^n_{2c}| \leq n^{10}\}.$$

$S$ first needs to visit $|k^n_{2a}|$ and $|k^n_{2c}|$ in order to visit both points $n^{10}$ and $-n^{10}$, when $|k^n_{2a}|, |k^n_{2c}| \leq n^{10}$ (since $S$ starts at the origin). Thus,

$$E^n_{k_a,c} \cap E^n_2 \subset E^n_{\text{visit}}.$$

Thus,

$$P(E^{nc}_{k_a,c}) + P(E^{nc}_2) \geq P(E^{nc}_{\text{visit}}).$$

If $P(E^{nc}_{k_a,c})$ and $P(E^{nc}_2)$ are both finitely summable over $n$, then $P(E^{nc}_{\text{visit}})$ is also. We prove that $P(E^{nc}_{k_a,c})$ is finitely summable and leave the proof that $P(E^{nc}_2)$ is finitely summable to the reader since it is very similar to the other one. Let $X^{R+}$, respectively, $X^{R-}$, be the first passage time of the random walk $R(k)_{k\in\mathbb{N}}$, respectively, $R(-k)_{k\in\mathbb{N}}$, at the point $3n$. Let $E^n_{R+} := \{X^{R+} \leq n^{10}\}$ and let $E^n_{R-} := \{X^{R-} \leq n^{10}\}$. Then, $E^n_{R+} \cup E^n_{R-} = E^n_{k_a,c}$. Thus, $P(E^{nc}_{R+}) + P(E^{nc}_{R-}) \geq P(E^{nc}_{k_a,c})$. By symmetry, $P(E^{nc}_{R+}) = P(E^{nc}_{R-})$. Thus, $2P(E^{nc}_{R+}) \geq P(E^{nc}_{k_a,c})$. Let $Z_i$ denote the first passage time of $\{R(k)\}_{k\in\mathbb{N}}$ at the point $i$. Let $X_i := Z_i - Z_{i-1}$. Then, $X^{R+} := \sum_1^{3n} X_i$ and $P(E^{nc}_{R+}) = P(\sum_1^{3n} X_i > n^{10}) \leq P((\sum_1^{3n} X_i)^{1/3} > n^3)$. For positive numbers $a_1, a_2, \dots, a_j$, we always have that $(\sum_{l=1}^j a_i)^3 \geq \sum_{l=1}^j (a_i)^3$. Thus $\sum_{i=1}^{3n} (X_i)^{1/3} \geq (\sum_{i=1}^{3n} X_i)^{1/3}$. It follows that $P(E^{nc}_{R+}) \leq P(\sum_{i=1}^{3n} (X_i)^{1/3} \geq n^3)$. By Chebyshev, we get

$$P(E^{nc}_{R+}) \leq \frac{3E[(X_1)^{1/3}]}{n^2}.$$

In [5] it is shown that $E[(X_i)^{1/3}]$ is finite. Thus, $P(E^{nc}_{R+})$ is finitely summable over $n$, which finishes this proof.

*High probability of $E^n_{\text{stopping}}$.* Let $E^n_3 := \{\forall\, i \leq \exp(n^{0.99}), t_{2ai} \leq \exp(n^{0.999})\}$. If up to time $t_{2ai}$ with $i = \exp(n^{0.99})$ the localization test makes no mistake in identifying exactly all the crossings by $R \circ S$ of $(0, 3n)$ which occur in the same place, then $E^n_{\text{stopping}}$ holds. Thus, if $t_{2ai} \leq \exp(n^{0.999})$ for $i = \exp(n^{0.99})$ and the localization test makes no mistake of this type up to time $\exp(n^{0.999})$, then $E^n_{\text{stopping}}$ holds. Let $E^n_{\text{test correct}}$ be the event that for all $z_a, z_b \in \mathbb{Z}$ with $0 < |z_a|, |z_b| \leq n^{0.999}$ and for all $0 < i_a, i_b \leq n^{0.999}$ the localization test makes no error when comparing the crossings $(t_{1a}, t_{2a})$ and $(t_{1b}, t_{2b})$. [Here $(t_{1a}, t_{2a})$ and $(t_{1b}, t_{2b})$ are defined as in Lemma 8: $(t_{1a}, t_{2a})$ is the $i_a$th crossing by $S$ of the $z_a$th crossing by $R$ of $(0, 3n)$ and $(t_{1b}, t_{2b})$ is the $i_b$th crossing by $S$ of the $z_b$th crossing by $R$ of $(0, 3n)$.] Up to time $\exp(n^{0.999})$, $S$ can cross a crossing by $R$ at most $\exp(n^{0.999})$ times. Thus, if $(t_{1a}, t_{2a})$ and $(t_{1b}, t_{2b})$ occur before time $\exp(n^{0.999})$, then $0 < i_a, i_b \leq n^{0.999}$. Furthermore, to reach the $z$th crossing $(k^n_{1z}, k^n_{2z})$, $S$ needs first to cross all the crossings $(k^n_{1z'}, k^n_{2z'})$ with $z'$ strictly between 0 and $z$. Thus up to time $\exp(n^{0.999})$ $S$ cannot reach any crossing $(k^n_{1z}, k^n_{2z})$ with $|z| > \exp(n^{0.999})$. If the crossings $(t_{1a}, t_{2a})$ and $(t_{1b}, t_{2b})$ occur before time $\exp(n^{0.999})$, we hence have that $0 < i_a, i_b \leq n^{0.999}$ and $0 < |z_a|, |z_b| \leq n^{0.999}$. Thus, $E^n_3$ and $E^n_{\text{test correct}}$ both hold; the localization test makes no mistake in identifying which of $(0, 3n)$ occur in the same place up to time $t_{2ai}$. In this case, $E^n_{\text{stopping}}$ holds. Thus,

$$E^n_3 \cap E^n_{\text{test correct}} \subset E^n_{\text{stopping}}.$$

It follows that

$$P(E^{nc}_3) + P(E^{nc}_{\text{test correct}}) \geq P(E^{nc}_{\text{stopping}}).$$

If $P(E^{nc}_3)$ and $P(E^{nc}_{\text{test correct}})$ are both finitely summable over $n$, then $P(E^{nc}_{\text{stopping}})$ is also. The proof that $P(E^{nc}_3)$ is finitely summable is very similar to the proof for $P(E^{nc}_{k_a,c})$, so we leave it to the reader. Let $E^n_{\text{test correct } i_a, i_b, z_a, z_b}$ be the event that the localization test recognizes correctly if with the crossings $(t_{1a}, t_{2a})$ and $(t_{1b}, t_{2b})$ we are in the $H_0$-case or not. By definition,

$$\bigcap E^n_{\text{test correct } i_a, i_b, z_a, z_b} = E^n_{\text{test correct}},$$

where the last intersection is taken over all $i_a, i_b, z_a, z_b$ such that $0 < |z_a|, |z_b| \leq n^{0.999}$ and $0 < i_a, i_b \leq n^{0.999}$. Thus,

$$\sum P(E^{nc}_{\text{test correct } i_a, i_b, z_a, z_b}) \geq P(E^{nc}_{\text{test correct}}),$$

where the sum is taken over the same domain as before the union. There are $n^{3.996}$ quadruples $(i_a, i_b, z_a, z_b)$ such that $0 < |z_a|, |z_b| \leq n^{0.999}$ and $0 < i_a, i_b \leq n^{0.999}$. By the large deviation principle and Lemma 8, the probability $P(E^{nc}_{\text{test correct } i_a, i_b, z_a, z_b})$ is exponentially small in $n$. Thus there exist $k > 0$ not depending on $n$ or on $(i_a, i_b, z_a, z_b)$ such that $P(E^{nc}_{\text{test correct } i_a, i_b, z_a, z_b}) \leq \exp(-kn)$. This implies that

$$P(E^{nc}_{\text{test correct}}) \leq n^{3.996} \cdot \exp(-kn).$$

Thus, $P(E^{nc}_{\text{test correct}})$ is finitely summable over $n$.

*High probability of $E^n_{\text{straight}}$.* Let $\bar{t}^n_{2ai}$ denote the 20,000th stopping time $t^n_{2ai}$. Thus, $\bar{t}^n_{2ai} := t^n_{2a(20{,}000 \cdot i)}$. Let $E^n_4$ be the event that there exists $i \leq n^{-20{,}000} \cdot e^{\bar{n}}$ and $s, r$ with $\bar{t}^n_{2ai} \leq s, r \leq \bar{t}^n_{2ai} + n^{220}$ such that $(s, r)$ is a straight crossing by $S$ of $(k^n_{1c}, k^n_{2a})$. We have that $E^n_4 \subset E^n_{\text{straight}}$. Let $E^n_5 := E^n_{k_a,c} \cap E^{\dot{n}}_{k_a,c}$. We find that the last inclusion implies

$$P(E^{nc}_4 \cap E^n_5) + P(E^{nc}_5) \geq P(E^{nc}_{\text{straight}}).$$

We already saw that $P(E^{nc}_5)$ is finitely summable over $n$. So it only remains to be proven that $P(E^{nc}_4 \cap E^n_5)$ is finitely summable over $n$. Let $X_i$ be the Bernoulli variable which is equal to 1 iff there exists $s, r$ with $\bar{t}^n_{2ai} \leq s, r \leq \bar{t}^n_{2ai} + n^{220}$ such that $(s, r)$ is a straight crossing by $S$ of $(k^n_{1c}, k^n_{2a})$. By the Markov property of the random walk $S$, we have that conditional under $\sigma(R(k)|k \in \mathbb{Z})$ the variables $X_i$ are i.i.d. Also, $E^n_5$ is $\sigma(R(k)|k \in \mathbb{Z})$-measurable. We are next going to evaluate the conditional probability: $P(X_1 = 1|R(k), k \in \mathbb{Z})$ when $E^n_5$ holds. When $E^n_5$ holds, then $|k^{\dot{n}}_{2a} - k^n_{2c}| \leq 2\dot{n}^{10}$. We have $2\dot{n}^{10} := 2n^{110}$. By definition at any time $t^{\dot{n}}_{2ai}$ the random walk $S$ is at the point $k^{\dot{n}}_{2a}$. By the local central limit theorem, when $|k^{\dot{n}}_{2a} - k^n_{2c}| \leq 2\dot{n}^{10}$, the probability that $S$ goes from $k^{\dot{n}}_{2a}$ to $k^n_{2c}$ in less than $\frac{1}{2}n^{220}$ steps is bigger than $k_2 \cdot n^{-110}$. (Here $k_2$ denotes a constant not depending on $n$ and not depending on $R$ as long as $R \in E^n_5$.) Crossing in a straight way to the point $k^n_{2a}$ right after the random walk $S$ is at the point $k^n_{2c}$, has probability bigger than $(\frac{1}{2})^{2n^{10}}$, when $|k^n_{2a} - k^n_{2c}| \leq 2n^{10}$. But, when $E^n_5$ holds, $|k^n_{2a} - k^n_{2c}| \leq 2n^{10}$. All this implies that when $E^n_5$ holds,

$$P(X_1 = 1|R(k), k \in \mathbb{Z}) \geq (k_2 n^{-110})(\tfrac{1}{2})^{2n^{10}}. \tag{16}$$

Let $\varepsilon_1 := (k_2 n^{-110})(\frac{1}{2})^{2n^{10}}$. Let $\hat{n} := n^{-20{,}000} \cdot e^{\bar{n}}$. Note that

$$E^{nc}_4 := \left\{ \sum_{i=1}^{\hat{n}} X_i = 0 \right\}.$$

Conditional under $\sigma(R(k)|k \in \mathbb{Z})$ the $X_i$'s are i.i.d. Thus,

$$P(E^{nc}_4|R(k), k \in \mathbb{Z}) = (1 - P(X_i = 1|R(k), k \in \mathbb{Z}))^{\hat{n}}.$$

Using (16), we get for $R \in E^n_5$,

$$P(E^{nc}_4|R(k), k \in \mathbb{Z}) \leq (1 - \varepsilon_1)^{\hat{n}}. \tag{17}$$

When $n$ goes to infinity, then $\varepsilon_1$ tends to zero. Thus, for $n$ big enough we get

$$(1 - \varepsilon_1)^{1/\varepsilon_1} \leq e^{-0.5}.$$

Applying this to (17) leads in the case that $R \in E^n_5$, to

$$P(E^{nc}_4|R(k), k \in \mathbb{Z}) \leq e^{-0.5\hat{n}\varepsilon_1}.$$

Integrating the last inequality over $E_5^n$ leads to

$$P(E_4^{nc} \cap E_5^n) \leq e^{-0.5\hat{n}\varepsilon_1}. \tag{18}$$

Recall that $\bar{n} := n^{10.89}$ and $\dot{n} := n^{11}$. In $\hat{n}$, the leading term is $e^{\bar{n}}$. In $\varepsilon_1$ the leading term is $e^{\ln(0.5)\cdot 2n^{10}}$. Since $n^{10.89} \gg n^{10}$ we get that $e^{\bar{n}} \gg e^{-\ln(0.5)\cdot 2n^{10}}$. This implies that the leading term in $\hat{n}\varepsilon_1$ is $e^{\bar{n}}$. Thus, the term on the right-hand side of (18) is finitely summable over $n$.

*High probability of $E^n_{\text{no other crossing by } R}$.* Let $n^* := n^{110} + n^{220}$. Let $(t^n_{111}, t^n_{211})$ designate the first crossing by $S$ of $(k^n_{11}, k^n_{21})$. Let $w^n_{11} := w_{(t^n_{111}, t^n_{211})}$. Define

$$E^n_{61} := \{\text{The only crossing } (k^n_{1z}, k^n_{2z}) \text{ with } 0 < |z| \leq n^* \\ \text{such that } w^R_{(k^n_{1z}, k^n_{2z})} \geq w^n_{11} \text{ is } (k^n_{11}, k^n_{21})\}.$$

Let $(t^n_{1(-1)1}, t^n_{2(-1)1})$ designate the first crossing by $S$ of $(k^n_{1(-1)}, k^n_{2(-1)})$. Let $w^n_{1(-1)} := w_{(t^n_{1(-1)1}, t^n_{2(-1)1})}$. Define

$$E^n_{6(-1)} := \{\text{The only crossing } (k^n_{1z}, k^n_{2z}) \text{ with } 0 < |z| \leq n^* \\ \text{such that } w^R_{(k^n_{1z}, k^n_{2z})} \geq w^n_{1(-1)} \text{ is } (k^n_{1(-1)}, k^n_{2(-1)})\}.$$

If $E^{\dot{n}}_{k_a,c}$ holds, then $|k^{\dot{n}}_{2a}| \leq n^{110}$. All the crossings $(k_1, k_2)$ concerned by the event $E^n_{\text{no other crossing by } R}$ are such that $|k_1 - k^{\dot{n}}_{2a}|, |k_2 - k^{\dot{n}}_{2a}| \leq n^{220}$. Thus, when $E^{\dot{n}}_{k_a,c}$ holds, then all the crossings concerned by $E^n_{\text{no other crossing by } R}$ are within $n^*$ of the origin. When we write those crossings in the form $(k^n_{1z}, k^n_{2z})$ they must be such that $|z| \leq n^*$. Thus, when $E^{\dot{n}}_{k_a,c}$ holds, the events $E^n_{61}$ and $E^n_{62}$ cover all the crossings involved in the definition of the event $E^n_{\text{no other crossing by } R}$. One of the crossings $(k^n_{1a}, k^n_{2a})$ and $(k^n_{1c}, k^n_{2c})$ is equal to $(k^n_{11}, k^n_{21})$ while the other one is equal to $(k^n_{1(-1)}, k^n_{2(-1)})$. Similarly, one of the crossings $(t^n_{1a}, t^n_{2a})$ and $(t^n_{1c}, t^n_{2c})$ is equal to $(t^n_{111}, t^n_{211})$ while the other one is equal to $(t^n_{1(-1)1}, t^n_{2(-1)1})$. Eventually, one of the words $w^n_a$ and $w^n_c$ is equal to $w^n_{11}$ while the other one is equal to $w^n_{1(-1)}$. This implies that when $E^{\dot{n}}_{k_a,c}$ holds, the events $E^n_{61}$ and $E^n_{62}$ jointly imply $E^n_{\text{no other crossing by } R}$. Thus,

$$E^n_{61} \cap E^n_{62} \cap E^{\dot{n}}_{k_a,c} \subset E^n_{\text{no other crossing by } R}.$$

It follows that

$$P(E^{nc}_{61}) + P(E^{nc}_{62}) + P(E^{\dot{n}c}_{k_a,c}) \geq P(E^{nc}_{\text{no other crossing by } R}).$$

We already saw that $P(E^{\dot{n}c}_{k_a,c})$ is finitely summable over $n$. By symmetry $P(E^{nc}_{61}) = P(E^{nc}_{62})$. Thus, it only remains to prove that $P(E^{nc}_{61})$ is finitely summable over $n$. Let

$$E^n_{61z} := \{w^R_{(k^n_{1z}, k^n_{2z})} \not\geq w^n_{11}\}.$$

We have

$$E_{61}^{n} := \bigcap_{0<|z|\leq n^*, z\neq 1} E_{61z}^{n}.$$

It follows that

$$P(E_{61}^{nc}) \leq \sum_{0<|z|\leq n^*, z\neq 1} P(E_{61z}^{nc}).$$

We saw in the proof of Lemma 8 the distribution, of $w^R_{(k_{1z}^n, k_{2z}^n)}$ does not depend on $z$. Thus, the expression on the right-hand side of the last inequality is equal to $(2n^* - 2)P(E_{612}^{nc})$. This yields

$$P(E_{61}^{nc}) \leq (2n^* - 2)P(E_{612}^{nc}). \tag{19}$$

We have that $E_{612}^{nc} = \{w^R_{(k_{12}^n, k_{22}^n)} \geq w_{11}^n\}$. Hence,

$$E_{612}^{nc} = \bigcap_{m=0}^{n-1} \{w^R_{(k_{12}^n, k_{22}^n)}(m) \geq w_{11}^n(m)\}.$$

As in the proof of Lemma 8, the bits of the word $w^R_{(k_{12}^n, k_{22}^n)}$ are i.i.d. as well as the bits of $w_{11}^n$ and $w^R_{(k_{12}^n, k_{22}^n)}$ is independent of $w_{11}^n$. This gives

$$P(E_{612}^{nc}) = \prod_{m=0}^{n-1} P(w^R_{(k_{12}^n, k_{22}^n)}(m) \geq w_{11}^n(m)) = P(w^R_{(k_{12}^n, k_{22}^n)}(1) \geq w_{11}^n(1))^n.$$

The probability $q := (w^R_{(k_{12}^n, k_{22}^n)}(1) \geq w_{11}^n(1))$ is strictly smaller than 1 and does not depend on $n$. Thus, the bound $(2n^* - 2)q^n$ on the left-hand side of (19) is finitely summable over $n$.

**5. Why the reconstruction of $\boldsymbol{\xi}$ works.** Our reconstruction algorithm constructs a scenery $\bar{\xi}$. The main result of this paper is that a.s. $\bar{\xi}$ is equivalent to $\xi$. This is also what we need to prove in this section. The reconstruction algorithm we propose constructs $\bar{\xi}$ by assembling (as explained in Section 3) the finite reconstructed pieces $\hat{\xi}^n$. The piece $\hat{\xi}^n$ is provided by the reconstruction algorithm at level $n$. The reconstruction algorithm at level $n$ tries to reconstruct the finite piece of the scenery $\xi$:

$$\xi^n := \xi(k_{1c}^n), \xi(k_{1c}^n + u), \xi(k_{1c}^n + 2u), \ldots, \xi(k_{1a}^n),$$

where $u := (k_{1a}^n - k_{1c}^n)/|k_{1a}^n - k_{1c}^n|$. We have proven in the last section that $(1 - P(\xi^n = \hat{\xi}^n))$ is finitely summable over $n$. It follows that a.s. $\xi^n = \hat{\xi}^n$ for all but a finite number of $n$'s. In Section 3 we have seen that the constructed scenery $\bar{\xi}$ is equivalent to $\xi$ as soon as for all but a finite number of $n$'s we have that $\xi^n \preccurlyeq_1 \xi^{n+1}$ and $\xi^n = \hat{\xi}^n$. It thus only remains to prove that a.s. for all but a finite number of $n$'s, $\xi^n \preccurlyeq_1 \xi^{n+1}$ holds. Define

$$\xi^n_{\text{inside}} := \xi(-n), \xi(-n+1), \xi(-n+2), \ldots, \xi(n)$$

and

$$\xi^n_{\text{outside}} := \xi(-n^{10}), \xi(-n^{10}+1), \xi(-n^{10}+2), \ldots, \xi(n^{10}).$$

By definition, $|k^n_{2a}|, |k^n_{2c}| \geq n$ from which it follows that $\xi^n_{\text{inside}} \preccurlyeq \xi^n$. On the other hand, if $E^n_{k_a,c}$ holds, then $|k^n_{2a}|, |k^n_{2c}| \leq n^{10}$ and $\xi^n \preccurlyeq \xi^n_{\text{outside}}$. Recall that $\xi^n \preccurlyeq \xi^{n+1}$ always holds by definition. Summing up: when $E^{n+1}_{k_a,c}$ holds, we find that

$$\xi^n_{\text{inside}} \preccurlyeq \xi^n \preccurlyeq \xi^{n+1} \preccurlyeq \xi^{n+1}_{\text{outside}}.$$

Next, note that if $\zeta_a, \zeta_b, \zeta_c, \zeta_d \in \bigcup_{l \in \mathbb{N}} \{0,1\}^l$ with $\zeta_a \preccurlyeq \zeta_b \preccurlyeq \zeta_c \preccurlyeq \zeta_d$ and $\zeta_a \preccurlyeq_1 \zeta_d$, then also $\zeta_b \preccurlyeq_1 \zeta_c$. Thus, when $E^{n+1}_{k_a,c}$ holds, if $\xi^n_{\text{inside}} \preccurlyeq_1 \xi^{n+1}_{\text{outside}}$, then also $\xi^n \preccurlyeq_1 \xi^{n+1}$. Let

$$E^n_{\text{unique}} := \{\xi^n_{\text{inside}} \preccurlyeq_1 \xi^{n+1}_{\text{outside}}\}.$$

We have shown that

$$E^n_{\text{unique}} \cap E^{n+1}_{k_a,c} \subset \{\xi^n \preccurlyeq_1 \xi^{n+1}\}.$$

For $(1 - P(\xi^n \preccurlyeq_1 \xi^{n+1}))$ to be finitely summable over $n$, it is thus enough that $P(E^{nc}_{\text{unique}})$ and $P(E^{c(n+1)}_{k_a,c})$ both are. We have already proven that the probability of the complement $P(E^{c(n+1)}_{k_a,c})$ is finitely summable over $n$. It remains to show that $P(E^{nc}_{\text{unique}})$ also is finitely summable. Let

$$E^n_{\text{unique, }+l} := \{\xi^n_{\text{inside}} \neq (\xi(l), \xi(l+1), \xi(l+2), \ldots, \xi(l+2n))\}$$

and

$$E^n_{\text{unique, }-l} := \{\xi^n_{\text{inside}} \neq (\xi(l), \xi(l-1), \xi(l-2), \ldots, \xi(l-2n))\}.$$

With this notation,

$$\bigcap_{l \neq -n, |l| \leq n^{10}} (E^n_{\text{unique, }+l} \cap E^n_{\text{unique, }-l}) \subset E^n_{\text{unique}}.$$

The last inclusion implies

$$\sum_{l \neq -n, |l| \leq n^{10}} P(E^{nc}_{\text{unique, }+l}) + P(E^{nc}_{\text{unique, }-l}) \geq P(E^{nc}_{\text{unique}}). \tag{20}$$

Because the scenery $\xi$ consists of i.i.d. Bernoulli variables with parameter $\frac{1}{2}$, we find that $P(E^{nc}_{\text{unique, }+l}) = P(E^{nc}_{\text{unique, }-l}) = (\frac{1}{2})^{2n}$. Furthermore, there are less than $2n^{10}$ elements in the set $\{l \neq -n, |l| \leq n^{10}\}$. This finishes the proof that the bound on the left-hand side of (20) is finitely summable over $n$.

**Acknowledgments.** I want to express my deepest gratitude to my thesis adviser Professor Harry Kesten, for providing such an interesting question and for all the time he spend helping me. I want to thank also Franz Merkl for always pushing me to finish writing up my articles in a readable way and for being a good friend. I thank the referee for motivating me to rewrite the paper in a more understandable way and for proofreading the new version carefully.

Department of Mathematics
University of Bielefeld
Postfach 100131
33501 Bielefeld
Germany
e-mail: matzing@mathematik.uni-bielefeld.de